\documentclass[a4paper]{article}

\usepackage{amsmath,amssymb}
\usepackage{a4wide}

\usepackage{tikz}
\usetikzlibrary{decorations.pathreplacing}
\usetikzlibrary{math}

\usepackage[titletoc, title]{appendix}

\usepackage{algorithm}
\usepackage{algpseudocode}
\usepackage[export]{adjustbox}
\usepackage{tabularx, booktabs}

\usepackage{graphicx}
\graphicspath{{./figures/}}

\newtheorem{remark1}{Remark}[section]

\newcommand{\argmin}[1]{\underset{#1}{\operatorname{arg}\!\operatorname{min}}\;}

\usepackage[normalem]{ulem} 

\newcounter{mainStokes}
\newcounter{subStokes}
\newcounter{tmpcount}

\begin{document}

\begin{center}

\begin{Large}
\textbf{A coupling concept for Stokes-Darcy systems: the ICDD method}
\end{Large}

\medskip

Marco Discacciati$^{1}$, Paola Gervasio$^{2}$

\medskip

${}^1$ Department of Mathematical Sciences, Loughborough University, Loughborough, United Kingdom, m.discacciati@lboro.ac.uk.

${}^2$ DICATAM, Universit\`a degli Studi di Brescia, Via Branze 38, 25123 Brescia, Italy, paola.gervasio@unibs.it.

\end{center}

\bigskip

\centerline{\textbf{Abstract}}

We present a coupling framework for Stokes-Darcy systems valid for arbitrary flow direction at low Reynolds numbers and for isotropic porous media. The proposed method is based on an overlapping domain decomposition concept to represent the transition region between the free-fluid and the porous-medium regimes. Matching conditions at the interfaces of the decomposition impose the continuity of velocity (on one interface) and pressure (on the other one) and the resulting algorithm can be easily implemented in a non-intrusive way. The numerical approximations of the fluid velocity and pressure obtained by the studied method converge to the corresponding counterparts computed by direct numerical simulation at the microscale, with convergence rates equal to suitable powers of the scale separation parameter $\varepsilon$ in agreement with classical results in homogenization.

\medskip

\emph{Keywords:} Stokes equations; Darcy equations; Domain decomposition methods; Multiscale analysis; Direct numerical simulations; Homogenization.

\section{Introduction}

Fluid flow systems formed by a free fluid region and an adjacent porous material through which the fluid can filtrate are of utmost relevance for a wide range of environmental, industrial and biomedical applications such as, e.g., membrane filtration for water purification, geophysical flows, neotissue growth in perfusion bioreactors, blood flow in biological tissues, water-gas flow in fuel cells, binary alloy solidification.
In general, due to the difficulty of reproducing the complex geometrical features of the porous medium at the microscale, mathematical models for filtration processes are not entirely based on the Navier-Stokes equations, but they rather include averaged macroscale models to represent the motion of the fluid inside the porous material.

A widely used upscaled modelling framework is the so-called two-domain approach, where the free-fluid region and the porous region are treated as two non-overlapping subdomains separated by an artificial sharp interface. Two different sets of governing equations are then used within each region, typically, the (Navier-)Stokes equations in the free fluid part and the Darcy equation \cite{Darcy:1856} or the Forchheimer equation \cite{Forchheimer:1901:ZVDI} in the porous medium. Specific coupling conditions must be imposed between the two sets of equations to correctly match the different fluid regimes. A classical two-domain model is defined by the Stokes and the Darcy equations coupled by the Beavers--Joseph--Saffman (BJS) conditions to impose mass conservation, balance of normal forces, and a slip condition on the interface \cite{Jager:1996:AnnPisa,Layton:2003:SINUM,Discacciati:2009:RMC}. These conditions, obtained experimentally by Beavers and Joseph \cite{Beavers:1967:JFM}, and simplified by Saffman \cite{Saffman:1971:SAM} to neglect the contribution of the smaller porous-medium tangential velocity at the interface, have also been given mathematical justification via homogenization theory and asymptotic analysis by J\"ager and Mikeli\'c in \cite{Jager:1996:AnnPisa,Jager:2000:SIAM}.  Although largely used in the literature, the BJS conditions were obtained in the case of a porous medium with quite large porosity and for cross-flow filtration, i.e., for a fluid flow parallel to the surface of the porous medium. Therefore, they could not be valid and are not justified in the case of arbitrary flow directions, for which finding suitable interface conditions is an active research area (see, e.g., \cite{Zampogna:2016:JFM, Lacis:2017:JFM, Naqvi:2021:IJMF,Eggenweiler:2021:MMS,Strohbeck:2023:TPM}). Models based on the BJS conditions are dependent on the position of the sharp interface (see,  e.g., \cite{Marciniak:2012:MMS,Eggenweiler:2022}) and they involve arbitrary coefficients that should be determined either by a detailed study of the physical characteristics of the porous medium or, within a  mathematical context, by solving a set of auxiliary ad-hoc problems across the interface. Alternative sets of interface conditions for systems of Stokes-Darcy type have also been obtained by methods of volume averaging and asymptotic expansion, e.g., in \cite{Angot:2017:PhysRevE,Angot:2021:AWR,Chandesris:2006:IJHMT,Chandesris:2009:TPM,OchoaTapia:1995:IJHMT}.

The sharp interface two-domain setting is obtained as the limiting case of collapsing to zero the width of the transition region that exists between the free fluid and the porous medium subregions, although physical quantities such as velocity and pressure undergo strong but still continuous variations inside such region \cite{LeBars:2006:JFM,Naqvi:2021:IJMF}.
As a result of this process achieved by volume averaging or homogenization techniques, the transition region itself is replaced by the sharp interface and the associated coupling conditions. Accurately accounting for the fluid variations within the transition region is crucial to ensure a reliable representation of the fluid flow. Indeed, despite the very small thickness of the transition region, typically of the size of the pores that characterize the porous material, the fluid behaviour therein significantly impacts both the free fluid and the porous regions. While in the two-domain approach, this is attempted by carefully selecting the position of the interface and by calibrating the coefficients appearing in the coupling conditions, other approaches are found in the literature.

The so-called one-domain method uses a single set of equations throughout the domain (typically, the Navier-Stokes-Brinkman equations \cite{Brinkman:1947:ASRA}) which rely on ad-hoc space-dependent physical parameters, like the effective fluid viscosity, the permeability and the porosity of the porous medium \cite{Chandesris:2006:IJHMT,Chandesris:2009:TPM,Cimolin:2013:ANM,Hernandez:2022:CES,Valdes:2013:AWR}, that must be carefully defined to model the transition layer. The characterization of these parameters is usually performed by matched asymptotic expansion methods, but their correct definition is far from being straightforward and this constitutes the principal limitation of the one-domain approach.

Recently, multi-domain methods have been proposed in \cite{Ruan:2023} and \cite{Kang:2024:IJMS}. In these works, the transition region is treated as an additional intermediate domain between the free fluid and the porous medium where the Brinkmann equations are used to model the fluid at the intermediate regime that occurs before the bulk region of the porous medium is reached. Coupling conditions at the two interfaces with the free fluid domain and the inner porous domain are then proposed. In \cite{Kang:2024:IJMS}, the first interface is placed on the upper surface of the porous medium, while the second one is set at depth equal to the square root of the permeability below the first one. This agrees with results in \cite{LeBars:2006:JFM}, where the authors also concluded that the Darcy regime in the porous medium holds only below a depth proportional to a few pore lengths. 

The modelling framework that we propose in this paper is also based on the idea of including a transition region of suitably small but positive thickness, but, differently from \cite{Ruan:2023,Kang:2024:IJMS}, we do not introduce any additional model therein. More precisely, we consider the case of laminar flow at low Reynolds numbers so that the Stokes equations and the Darcy equations can be used in the free fluid region and in the porous medium, respectively. Then, following \cite{LeBars:2006:JFM}, we consider the Stokes equations to be valid at leading order inside the whole transition layer that, in our framework, forms an overlapping region between the free fluid domain and the porous medium domain. Simple leading-order matching conditions imposing the continuity of the velocity on the upper interface and of the pressure on the lower one are set to complete the definition of the model. The resulting coupled problem does not rely on any assumption concerning the flow direction, and it does not require the solution of auxiliary ad-hoc problems to identify coupling parameters. Despite its simple setting, we show that, from a physical point of view, the proposed coupling concept provides a leading-order approximation of the Stokes-Darcy velocity and pressure that respectively converge to the microscale velocity and pressure computed by direct numerical simulation. The estimated orders of convergence agree with classical theoretical results from homogenization theory as proved in, e.g., \cite{Marciniak:2012:MMS} for the cross-flow filtration case. 

From a computational point of view, our coupling framework leads to a robust overlapping domain decomposition method, the so-called Interface Control Domain Decomposition (ICDD) method that was first introduced in \cite{Discacciati:2013:SICON} and then analyzed for the Stokes-Darcy problem in \cite{Discacciati:2016:SINUM}. The algorithm can be formulated in a completely non-intrusive way and it can be easily implemented using existing CFD software.

The paper is organised as follows. Section \ref{sec:stokesDarcy} introduces the ICDD approach for the Stokes-Darcy problem and it details its computer implementation. Section \ref{sec:microscale} focusses on the physical interpretation of the overlapping approach and it gives a practical strategy to define the interfaces between the fluid region and the porous medium domain. Finally, Section \ref{sec:validation} provides numerical evidence of the effectiveness and accuracy of the ICDD method by comparing the numerical results computed by this strategy with those obtained by direct numerical simulations at the microscale. Error estimates in line with existing literature in homogenization theory are finally presented and discussed.

\section{Macroscale Stokes-Darcy model with overlap: the ICDD approach}\label{sec:stokesDarcy}

\begin{figure}[bht]
\begin{center}
 	\begin{tikzpicture}[scale=0.75]
        
        \fill[gray!20] (-0.5,2.2) rectangle (4.5,2.6);
        \fill[gray!40] (-0.5,-0.5) rectangle (4.5,2.2);
        
        \draw (-0.5,-0.5)--(4.5,-0.5);
        \draw (-0.5,-0.5)--(-0.5,5.5);
        \draw (-0.5,5.5)--(4.5,5.5);
        \draw (4.5,-0.5)--(4.5,5.5);
        \draw[dashed] (-0.5,2.2)--(4.5,2.2); \node[black] at (3.5, 1.8) {${\Gamma}_f$};
        \draw[dashed] (-0.5,2.6)--(4.5,2.6); \node[black] at (3.5, 2.9) {${\Gamma}_p$};

		\node[black] at (0.3, 1) {${\Omega}_p$};
		\draw[<->, black,dotted]  (0.7, -0.5) -- (0.7, 2.6);
		
		\node[black] at (1.2, 4) {${\Omega}_f$};
 		\draw[<->, black,dotted]  (1.6, 2.2) -- (1.6, 5.5);

		\node[black] at (2.0, -0.9) {$\Omega$};
 	\end{tikzpicture}
	\end{center}
  \caption{Schematic representation of the macroscale domain $\Omega$.}
  \label{fig:domainMicroMacro}
\end{figure}
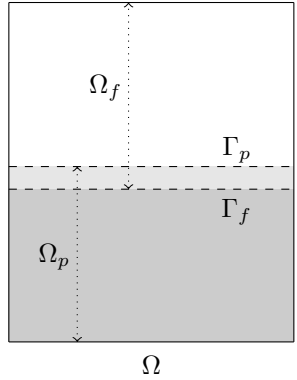

In this section, we formulate a macroscale model by considering the Stokes equations and the Darcy equations in two overlapping regions and by imposing suitable matching conditions between the fluid velocity and pressure at the boundary of the overlapping domain.
To this aim, let $\Omega$ be the computational domain split into two overlapping subdomains $\Omega_f$ and $\Omega_p$ as shown in Fig.~\ref{fig:domainMicroMacro}. The internal boundaries of $\Omega_f$ and $\Omega_p$ are denoted by $\Gamma_f$ and $\Gamma_p$ and they are called interfaces.

\smallskip

We assume that $\Omega$ is filled by an incompressible fluid characterized by low Reynolds number so that, in the fluid domain $\Omega_f$, the fluid behaviour can be modelled by the Stokes equations
\setcounter{mainStokes}{\value{equation}}
\begin{subequations}\label{eq:stokesFluid}
  \begin{eqnarray}
      -\mu \, \Delta \mathbf{u}_f + \nabla p_f = \mathbf{f}   && \text{in } \Omega_f \, ,\label{eq:stokesFluid_momentum}\\
                     \nabla \cdot \mathbf{u}_f = 0            && \text{in } \Omega_f \, , \\
                                  \mathbf{u}_f = \mathbf{u}_D && \text{on } \partial\Omega_f \setminus \Gamma_f\,, \label{eq:stokesFluid_bc}
  \end{eqnarray}
  \setcounter{subStokes}{\value{equation}}
\end{subequations}
where $\mathbf{u}_f$ and $p_f$ are the fluid velocity and pressure, respectively, $\mu$ is the fluid dynamic viscosity,    
$\textbf{f}=\rho\textbf{b}$ is an external force per unit volume, 
with $\rho$ the density of the fluid 
and $\textbf{b}$ a suitable acceleration field, 
and $\mathbf{u}_D$ is an assigned boundary velocity field.

In the porous medium domain $\Omega_p$, the fluid is described by Darcy's equations
\begin{subequations}\label{eq:darcy}
  \begin{eqnarray}
  \mu\,\mathbf{u}_p = - \mathbf{K} (\nabla p_p -\mathbf{f}) && \text{in } \Omega_p \, , \label{eq:darcy1} \\
                \nabla \cdot \mathbf{u}_p = 0 && \text{in } \Omega_p \, ,\\
             \mathbf{u}_p\cdot \mathbf{n} = 0 && \text{on } \partial\Omega_p \setminus \Gamma_p\,, \label{eq:darcy_bc}
  \end{eqnarray}
\end{subequations}
where $\mathbf{u}_p$ and $p_p$ are the Darcy's velocity and pressure, $\mathbf{K}$ is the permeability tensor 
and $\mathbf{n}$ is the unit normal vector on $\partial\Omega_p$ directed outwards of $\Omega_p$.

\smallskip

To complete the definition of the macroscopic model, we need to identify a suitable location of the interfaces $\Gamma_f$ and $\Gamma_p$ and the matching conditions between the Stokes and the Darcy variables to impose therein.

\smallskip

The position of the interface $\Gamma_p$ coincides with the top surface that delimits the microscale solid inclusions that form the porous material in $\Omega_p$ (see Fig.~\ref{fig:domainMicroscopic}, Sect.~\ref{sec:microscale}). 
The matching condition on $\Gamma_p$ plays the role of a boundary condition for Darcy's problem \eqref{eq:darcy}. Considering that Darcy's regime is mainly driven by the fluid pressure and its gradient, we set a condition on Darcy's pressure $p_p$ that imposes that this is equal to the Stokes pressure:
\begin{equation}\label{eq:continuityPressure}
    p_p = p_f \quad \mbox{on } \Gamma_p\, ,
\end{equation}
under the reasonable assumption that the pressures $p_f$ and $p_p$ remain bounded in the neighbourhood of the interfaces $\Gamma_f$ and $\Gamma_p$.

In our modelling framework, the overlapping region $\Omega_f \cap \Omega_p$ is used to account for the well-known presence of a transitional regime between a free flow and a porous medium, as we will discuss in details in Sect.~\ref{sec:microscale}. Therefore, the interface $\Gamma_f$ must be accurately placed below $\Gamma_p$ in order to capture the reduction of the magnitude of the Stokes velocity occurring when the fluid enters the porous material, but it must not be too low inside $\Omega_p$ to avoid reaching the bulk region of the porous medium where Darcy's regime holds and the velocity magnitude is very small. A practical way to correctly define $\Gamma_f$ will be discussed in Sect.~\ref{sec:microscale}.

The matching condition on $\Gamma_f$ defines the boundary condition for the Stokes problem \eqref{eq:stokesFluid}, and we choose to impose the continuity of the whole velocity field, i.e.,
\begin{equation}\label{eq:continuityVelocity}
    \mathbf{u}_f = \mathbf{u}_p \quad \mbox{on } \Gamma_f\, .
\end{equation}

Since we are imposing Dirichlet conditions on the whole boundary of the subdomain $\Omega_f$, to guarantee the uniqueness of the pressure $p_f$ we complete problem \eqref{eq:stokesFluid} with the null--average condition
\setcounter{tmpcount}{\value{equation}} 
\setcounter{equation}{\value{mainStokes}}
\begin{subequations}
\setcounter{equation}{\value{subStokes}}
\begin{equation}\label{eq:nullPressureMedia}
   \int_{\Omega_f} p_f =0.
\end{equation}
\setcounter{subStokes}{\value{equation}}     
\end{subequations}

\setcounter{equation}{\value{tmpcount}} 

The coupling conditions \eqref{eq:continuityPressure} and \eqref{eq:continuityVelocity} guarantee the well-posedness of the local Stokes and Darcy problems \eqref{eq:stokesFluid} and \eqref{eq:darcy} as well as of the global macroscopic model \eqref{eq:stokesFluid}--\eqref{eq:continuityVelocity} in $\Omega$. Indeed, this macroscopic model was initially proposed in \cite{Discacciati:2016:SINUM} where it was shown that it is equivalent to a well-posed optimal control problem with interface controls corresponding to the trace of the Stokes velocity $\mathbf{u}_f$ on $\Gamma_f$ and to the trace of the Darcy pressure $p_p$ on $\Gamma_p$, and with an interface cost functional that imposes the matching conditions \eqref{eq:continuityPressure} and \eqref{eq:continuityVelocity} weakly in a suitable trace space. 

From a computational point of view, the model defines the Interface Control Domain Decomposition (ICDD) method, an overlapping domain decomposition algorithm that requires solving iteratively and independently the Stokes problem \eqref{eq:stokesFluid} (completed with the boundary condition \eqref{eq:continuityVelocity}) and the Darcy problem \eqref{eq:darcy} (completed with the boundary condition \eqref{eq:continuityPressure}) (see Sect. \ref{sec:ICDD} for an efficient solution algorithm). At convergence, we define the global fluid velocity and pressure in $\Omega$ as
\begin{equation}\label{eq:solICDD}
  \mathbf{u} =
   \left\{
     \begin{array}{ll}
      \mathbf{u}_f & \text{ in } \Omega_f \\
      \mathbf{u}_p & \text{ in } \Omega_p \setminus (\Omega_f \cap \Omega_p)
     \end{array}
   \right.
   \quad \text{and} \quad
   p =
   \left\{
     \begin{array}{ll}
      p_f & \text{ in } \Omega_f \\
      p_p & \text{ in } \Omega_p \setminus (\Omega_f \cap \Omega_p) \, .
     \end{array}
   \right.
\end{equation}
In \eqref{eq:solICDD}, the Stokes velocity $\mathbf{u}_f$ and pressure $p_f$ are retained in the overlapping region instead of their Darcy counterparts because the overlapping region is purposefully designed so that the Stokes regime remains valid therein.

In the rest of the paper, we will refer to the macroscale Stokes-Darcy model with overlap \eqref{eq:stokesFluid}--\eqref{eq:continuityVelocity} as the \emph{ICDD model.} A thorough physical justification of the ICDD model will be provided in Sect.~\ref{sec:microscale}, while in the rest of this section we focus on practical aspects of its implementation.

\subsection{Practical implementation of the ICDD method for Stokes-Darcy}\label{sec:ICDD}

The Interface Control Domain Decomposition (ICDD) method \cite{Discacciati:2013:SICON,Discacciati:2013:JCSMD,Discacciati:2014:IJNMF,Discacciati:2016:SINUM} was designed to solve boundary value problems governed by partial differential equations using overlapping domain decompositions of the domain of interest. The method is based on an optimal control problem that imposes the continuity of selected physical variables at the interfaces of the overlapping region, and it can be applied to problems characterised by single or multiple physics. In this section, we explain how to practically implement the ICDD method to solve the Stokes--Darcy model with overlap \eqref{eq:stokesFluid}--\eqref{eq:continuityVelocity} and we refer to \cite{Discacciati:2016:SINUM} for its theoretical analysis.

Consider two computational meshes, one in $\Omega_f$ and one in $\Omega_p$ and, for simplicity of exposition, let us assume that the two meshes are conformal on the closed overlapping region $\overline{\Omega_f\cap\Omega_p}$, so that the nodes on $\Gamma_f$ (respectively, $\Gamma_p$) belong to the mesh defined in $\overline \Omega_p$ (respectively, $\overline \Omega_f$). 
To guarantee that ICDD solves correctly the problem (\ref{eq:stokesFluid})--(\ref{eq:continuityVelocity}), we require that $\overline{\Gamma}_f$ and $\overline{\Gamma}_p$ are disjoint (see \cite{Discacciati:2016:SINUM}). The Stokes and the Darcy equations are discretized by a suitable mesh-based method such as the Galerkin Finite Element Method (FEM) or the Spectral Element Methods (SEM).  In this work, we have used stabilized SEM both for Stokes \cite{Gervasio:1998:NMPDE} and Darcy equations, in the latter case by adapting the stabilized FEM method \cite{Masud:2002:CMAME} to SEM. (We refer, e.g., to \cite{Quarteroni:1994, Brooks:1982:CMAME, Franca:1992:CMAME} for other stabilization methods for the FEM discretization of Stokes equations.) The Darcy equations \eqref{eq:darcy} could be reformulated as a second-order elliptic equation for the pressure $p_p$, but since the velocity $\mathbf{u}_p$ on $\Gamma_p$ is needed at each iteration of the ICDD method, we have chosen to directly use the mixed formulation \eqref{eq:darcy} instead of reconstructing the velocity field at each iteration.

At the algebraic level, the Degrees of Freedom (DoFs) associated with the nodes on the interface $\Gamma_f$ (respectively, $\Gamma_p$) are separated from those associated with the other nodes in $\overline\Omega_f\setminus\Gamma_f$ (respectively, $\overline\Omega_p\setminus\Gamma_p$).
For the Stokes problem, let $\textsf{u}_f$ and $\textsf{p}_f$ denote the arrays of the DoFs of the discrete velocity and pressure associated with the mesh nodes in $\overline\Omega_f\setminus\Gamma_f$, and let $\textsf{g}_f$ be the array of the DoFs of the discrete velocity field associated with the nodes on $\Gamma_f$.
Similarly for the Darcy problem, $\textsf{u}_p$ and $\textsf{p}_p$ are the arrays of the discrete velocity and pressure DoFs associated with the mesh nodes in $\overline\Omega_p\setminus\Gamma_p$, and $\textsf{g}_p$ is the array of the DoFs of the discrete pressure at the nodes on $\Gamma_p$.

Using standard FEM notation (see, e.g., \cite{Quarteroni:1994}), the algebraic counterpart of problem \eqref{eq:stokesFluid}--\eqref{eq:continuityVelocity} reads:
\begin{equation}\label{eq:algebraic-stokes-darcy}
\left[\begin{array}{cccc|cc}
 \textsf{A}_f & \textsf{B}_f^t &  \textsf{0} & \textsf{0}  & \textsf{A}_{\Gamma_f} & \textsf{0} \\[0.8mm]
\textsf{B}_f & \textsf{C}_f & \textsf{0} & \textsf{0} & \textsf{0} & \textsf{0} \\[0.8mm]
\textsf{0} & \textsf{0} & \textsf{A}_p & \textsf{B}_p^t & \textsf{0} &  \textsf{A}_{\Gamma_p} \\[0.8mm]
\textsf{0} & \textsf{0} &\textsf{B}_p&  \textsf{C}_p & \textsf{0} & \textsf{0} \\[0.8mm]
\hline
& & & & &\\[-2mm]
\textsf{0} & \textsf{0} & -\textsf{R}_{\Gamma_f p} & \textsf{0} & \textsf{I}_f & \textsf{0}\\[0.8mm]
\textsf{0} & -\textsf{R}_{\Gamma_p f} & \textsf{0} & \textsf{0} & \textsf{0} & \textsf{I}_p\\[0.8mm]
\end{array}\right]
\left[\begin{array}{l}
\textsf{u}_f \\[0.8mm]
\textsf{p}_f \\[0.8mm]
\textsf{u}_p \\[0.8mm]
\textsf{p}_p \\[0.8mm]
\hline
\\[-2mm]
\textsf{g}_f\\[0.8mm]
\textsf{g}_p\\[0.8mm]
\end{array}\right]=
\left[\begin{array}{c}
\textsf{f}_f \\[0.8mm]
 \textsf{0} \\[0.8mm]
\textsf{f}_p \\[0.8mm]
\textsf{0} \\[0.8mm]
\hline
\\[-2mm]
\textsf{0}\\[0.8mm]
\textsf{0}\\[0.8mm]
\end{array}
\right] \, .
\end{equation}

For $\star\in\{f,p\}$, $\textsf{A}_\star$ is the square stiffness matrix associated with the DoFs in $\overline\Omega_\star\setminus\Gamma_\star$, while $\textsf{A}_{\Gamma_\star}$ is the rectangular block of the stiffness matrix whose rows are associated with the DoFs in $\overline\Omega_\star\setminus\Gamma_\star$ and whose columns are associated with the DoFs on $\Gamma_\star$. The matrix $\textsf{B}_\star$ corresponds to the discretization of the divergence-free term in $\Omega_\star$, $\textsf{C}_\star$ could be null depending both on the discretization and the stabilization adopted, while $\textsf{I}_\star$ is the identity matrix whose size is equal to the number of DoFs on $\Gamma_\star$.
Moreover, $\textsf{R}_{\Gamma_p f}$ (respectively, $\textsf{R}_{\Gamma_f p}$) is a rectangular restriction matrix with entries either 0 or 1 that, for every array of nodal values in $\overline\Omega_f\setminus\Gamma_f$ (respectively, $\overline\Omega_p\setminus\Gamma_p$), returns the array of nodal values on $\Gamma_p$ (respectively, $\Gamma_f$). These restriction matrices are well defined because $\overline{\Gamma}_f$ and $\overline{\Gamma}_p$ are disjoint. Finally, the right-hand side arrays $\textsf f_f$ and $\textsf f_p$ take into account both the external force $\mathbf{f}$ and the boundary conditions \eqref{eq:stokesFluid_bc} and \eqref{eq:darcy_bc}. Remark that the last two equations in \eqref{eq:algebraic-stokes-darcy} respectively impose the matching conditions \eqref{eq:continuityVelocity} and \eqref{eq:continuityPressure} at the algebraic level.

When the thickness of the overlapping region $\Omega_f \cap \Omega_p$ is very tiny, solving system \eqref{eq:algebraic-stokes-darcy} could become quite cumbersome, indeed, the number of iterations required to solve this linear system by an iterative method can grow when the overlap thickness decreases. To overcome this drawback, we reformulate the coupled system \eqref{eq:stokesFluid}--\eqref{eq:continuityVelocity} in an equivalent form as follows (see \cite{Discacciati:2016:SINUM}). Let $\mathbf{g}_f$ denote the unknown trace of the velocity $\mathbf{u}_f$ on $\Gamma_f$ and $g_p$ the unknown trace of the pressure $p_p$ on $\Gamma_p$. Then, we consider the system
\begin{eqnarray}\label{eq:OS-icdd}
\begin{array}{ll}
\mbox{\emph{Primal problems:}}\\[2mm]
\left\{\begin{array}{rl}
-\mu\,\Delta \mathbf{u}_f+\nabla p_f=\mathbf{f} & \mbox{ in }\Omega_f,\\
\nabla\cdot \mathbf{u}_f=0 & \mbox{ in }\Omega_f,\\
\mathbf{u}_f = \mathbf{g}_f & \mbox{ on }\Gamma_f,\\
\mathbf{u}_f=\mathbf{u}_D & \mbox{ on }\partial\Omega_f\setminus\Gamma_f,\\
\int_{\Omega_f}p_f=0,
\end{array}\right.
&
\left\{\begin{array}{rl}
\mu\,\mathbf{u}_p=-\mathbf{K}(\nabla p_p-\mathbf{f}) & \mbox{ in }\Omega_p,\\
\nabla\cdot \mathbf{u}_p=0 & \mbox{ in }\Omega_p,\\
p_p= g_p & \mbox{ on }\Gamma_p,\\
\mathbf{u}_p\cdot\mathbf{n}=0 & \mbox{ on }\partial\Omega_p\setminus\Gamma_p;
\end{array}\right.
\\
\\ 
\mbox{\emph{Dual problems:}}\\[2mm]
\left\{\begin{array}{rl}
-\mu\,\Delta \mathbf{w}_f+\nabla q_f=\mathbf{0} & \mbox{ in }\Omega_f,\\
\nabla\cdot \mathbf{w}_f=0 & \mbox{ in }\Omega_f,\\
\mathbf{w}_f = \mathbf{g}_f-\mathbf{u}_p & \mbox{ on }\Gamma_f,\\
\mathbf{w}_f=\mathbf{0} & \mbox{ on }\partial\Omega_f\setminus\Gamma_f,\\
\int_{\Omega_f}q_f=0,
\end{array}\right.
&
\left\{\begin{array}{rl}
\mu\,\mathbf{w}_p=-\mathbf{K}\nabla q_p & \mbox{ in }\Omega_p,\\
\nabla\cdot \mathbf{w}_p=0 & \mbox{ in }\Omega_p,\\
q_p= g_p-p_f & \mbox{ on }\Gamma_p,\\
\mathbf{w}_p\cdot\mathbf{n}=0 & \mbox{ on }\partial\Omega_p\setminus\Gamma_p;
\end{array}\right.
\\
\\ 
\mbox{\emph{Interface equations:}}\\[2mm]
\quad \mathbf{g}_f-\mathbf{u}_p+\mathbf{w}_p =0 \quad  \mbox{ on }\Gamma_f & 
\quad q_f+g_p-p_f =0 \quad \mbox{ on }\Gamma_p.

\end{array}
\end{eqnarray}

System (\ref{eq:OS-icdd}) is the optimality system associated with an optimal control problem aiming at minimizing the jump between the velocities on $\Gamma_f$ and the jump between the pressures on $\Gamma_p$. (If the boundary conditions on $\partial\Omega_\star\setminus\Gamma_\star$ (with $\star\in\{f,p\}$) for the primal problems are replaced by other boundary conditions, e.g., Neumann, the corresponding boundary conditions in the dual problems must be of the same type but homogeneous, with possible removal of the null average condition for the pressure.)
This system is solved iteratively (as we will show hereafter) with the dual problems playing the role of preconditioners stabilizing the iterative process, especially when the overlap thickness is very tiny. 
Notice that, at convergence, the solutions of the dual problems are null, $\mathbf{g}_f={\mathbf{u}_p}_{|\Gamma_f}$, and 
$g_p={p_f}_{|\Gamma_p}$, thus we recover the equivalence between (\ref{eq:stokesFluid})--(\ref{eq:continuityVelocity}) and (\ref{eq:OS-icdd}). For more details, we refer to \cite{Discacciati:2016:SINUM}.

Using the notations introduced before, the algebraic counterpart of system (\ref{eq:OS-icdd}) reads:

\begin{equation}\label{eq:algebraic-OS-icdd}
\hspace*{-3mm}
\left[\begin{array}{cccc|cccc|cc}
\textsf{A}_f & \textsf{B}_f^t &  \textsf{0} & \textsf{0}   &
\textsf{0} & \textsf{0}& \textsf{0} & \textsf{0} &
 \textsf{A}_{\Gamma_f} & \textsf{0} \\[0.8mm]
\textsf{B}_f & \textsf{C}_f & \textsf{0} & \textsf{0}  &
\textsf{0} & \textsf{0}& \textsf{0} & \textsf{0} &
 \textsf{0} & \textsf{0} \\[0.8mm]
\textsf{0} & \textsf{0} & \textsf{A}_p & \textsf{B}_p^t  &
\textsf{0} & \textsf{0}& \textsf{0} & \textsf{0} &
 \textsf{0} &  \textsf{A}_{\Gamma_p} \\[0.8mm]
\textsf{0} & \textsf{0} &\textsf{B}_p&  \textsf{C}_p  &
\textsf{0} & \textsf{0}& \textsf{0} & \textsf{0} &
 \textsf{0} & \textsf{0} \\[0.8mm]
\hline
& & & & & & & & &\\[-2mm]
\textsf{0} & \textsf{0} & -\textsf{A}_{\Gamma_f}\textsf{R}_{\Gamma_f p} &
\textsf{0}  &
\textsf{A}_f & \textsf{B}_f^t &  \textsf{0} & \textsf{0} &
\textsf{A}_{\Gamma_f} & \textsf{0} \\[0.8mm]
\textsf{0} & \textsf{0}& \textsf{0} & \textsf{0} &
\textsf{B}_f & \textsf{C}_f & \textsf{0} & \textsf{0}  &
 \textsf{0} & \textsf{0} \\[0.8mm]
\textsf{0} & -\textsf{A}_{\Gamma_p}\textsf{R}_{\Gamma_p f} & \textsf{0} & \textsf{0} &
\textsf{0} & \textsf{0} &\textsf{A}_p & \textsf{B}_p^t  &
 \textsf{0} &  \textsf{A}_{\Gamma_p} \\[0.8mm]
\textsf{0} & \textsf{0}& \textsf{0} & \textsf{0} &
\textsf{0} & \textsf{0} &\textsf{B}_p&  \textsf{C}_p  &
 \textsf{0} & \textsf{0} \\[0.8mm]
\hline
& & & & & & & & &\\[-2mm]
\textsf{0} & \textsf{0} & -\textsf{R}_{\Gamma_f p} &
\textsf{0}  &
\textsf{0} & \textsf{0}&\textsf{R}_{\Gamma_f p}& \textsf{0}  &
\textsf{I}_f& \textsf{0}\\[0.8mm]
\textsf{0} &-\textsf{R}_{\Gamma_p f} &\textsf{0} & \textsf{0} &
\textsf{0} &\textsf{R}_{\Gamma_p f} &\textsf{0} & \textsf{0}&
\textsf{0} &\textsf{I}_p\\[0.8mm]
\end{array}\right]
\left[\begin{array}{l}
\textsf{u}_f \\[0.8mm]
\textsf{p}_f \\[0.8mm]
\textsf{u}_p \\[0.8mm]
\textsf{p}_p \\[0.8mm]
\hline
\\[-2mm]
\textsf{w}_f \\[0.8mm]
\textsf{q}_f \\[0.8mm]
\textsf{w}_p \\[0.8mm]
\textsf{q}_p \\[0.8mm]
\hline
\\[-2mm]
\textsf{g}_f\\[0.8mm]
\textsf{g}_p\\[0.8mm]
\end{array}\right]=
\left[\begin{array}{c}
\textsf{f}_f \\[0.8mm]
 \textsf{0} \\[0.8mm]
\textsf{f}_p \\[0.8mm]
\textsf{0} \\[0.8mm]
\hline
\\[-2mm]
\textsf{0} \\[0.8mm]
\textsf{0} \\[0.8mm]
\textsf{0} \\[0.8mm]
\textsf{0} \\[0.8mm]
\hline
\\[-2mm]
\textsf{0}\\[0.8mm]
\textsf{0}\\[0.8mm]
\end{array}
\right]
\end{equation}

Notice that the matrix blocks $(1,1)$ and $(2,2)$ in \eqref{eq:algebraic-OS-icdd} are identical because the dual problems are of the same nature as the primal ones, hence no auxiliary matrices must be built and stored for the dual problems.
After introducing the matrices
\begin{eqnarray*}\label{eq:matricesRM}
\textsf{R}_\Gamma=
\left[\begin{array}{cccc}
\textsf{0} & \textsf{0} & \textsf{R}_{\Gamma_f p} & \textsf{0}\\
\textsf{0} & \textsf{R}_{\Gamma_p f} & \textsf{0} & \textsf{0}
\end{array}\right], 
\qquad
\textsf{I}_\Gamma=
\left[\begin{array}{cc}
\textsf{I}_{f} & \textsf{0}\\
\textsf{0} & \textsf{I}_{p}
\end{array}\right],
\end{eqnarray*}
\begin{eqnarray*}\label{eq:matricesA}
\textsf{A}=
\left[\begin{array}{cccc}
\textsf{A}_f & \textsf{B}_f^t & \textsf{0} & \textsf{0}\\
\textsf{B}_f & \textsf{C}_f   & \textsf{0} & \textsf{0}\\
\textsf{0} & \textsf{0} & \textsf{A}_p & \textsf{B}_p^t \\
\textsf{0} & \textsf{0} & \textsf{B}_p & \textsf{C}_p
\end{array}
\right],
\qquad
\textsf{A}_\Gamma=
\left[\begin{array}{cc}
\textsf{A}_{\Gamma_f} & \textsf{0}\\
\textsf{0} & \textsf{0}\\
\textsf{0} & \textsf{A}_{\Gamma_p}\\
\textsf{0} & \textsf{0}
\end{array}\right],
\end{eqnarray*}
and the arrays
\begin{eqnarray*}\label{eq:arraysUWg}
\begin{array}{cccc}
\textsf{U}=\left[\begin{array}{c}
\textsf{u}_f\\
\textsf{p}_f\\
\textsf{u}_p\\
\textsf{p}_p
\end{array}\right],
&
\textsf{W}=\left[\begin{array}{c}
\textsf{w}_f\\
\textsf{q}_f\\
\textsf{w}_p\\
\textsf{q}_p
\end{array}\right],
&
\textsf{F}=\left[\begin{array}{c}
\textsf{f}_f\\
\textsf{0}\\
\textsf{f}_p\\
\textsf{0}
\end{array}\right],
&
\textsf{g}=\left[\begin{array}{c}
\textsf{g}_f\\
\textsf{g}_p
\end{array}\right],
\end{array}
\end{eqnarray*}
system (\ref{eq:algebraic-OS-icdd}) can be written in the compact form
\begin{eqnarray}\label{eq:algebraic-OS-icdd-compact}
\left[\begin{array}{ccc}
\textsf{A} & \textsf{0}  & \textsf{A}_\Gamma\\
-\textsf{A}_\Gamma \textsf{R}_\Gamma & \textsf{A} & \textsf{A}_\Gamma \\
-\textsf{R}_\Gamma & \textsf{R}_\Gamma & \textsf{I}_\Gamma
\end{array}\right]
\left[\begin{array}{c}
\textsf{U}\\
\textsf{W}\\
\textsf{g}
\end{array}\right]
=
\left[\begin{array}{c}
\textsf{F}\\
\textsf{0}\\
\textsf{0}
\end{array}\right].
\end{eqnarray}

Instead of solving the global system \eqref{eq:algebraic-OS-icdd-compact}, we solve iteratively its Schur-complement system 
\begin{equation}\label{eq:Schur}
\textsf{S}\,\textsf{g}=\textsf{b}
\end{equation}
built with respect to the interface DoFs array $\textsf{g}$. This approach makes the implementation of the ICDD method completely non-intrusive and suitable for reusing existing computational codes specifically developed for the Stokes and Darcy problems. Indeed, the Schur-complement matrix $\textsf{S}$ and right-hand side $\textsf{b}$ need not be assembled, but the system (\ref{eq:Schur}) is tackled by solving iteratively and separately the local Stokes and Darcy subproblems and by exchanging information across the interfaces. More precisely, there holds
\begin{equation}\label{eq:S}
\textsf{S}=\textsf{I}_\Gamma+\textsf{R}_\Gamma
\left[\begin{array}{cc} \textsf{I}_\Gamma & -\textsf{I}_\Gamma\end{array}\right]
\left[\begin{array}{cc}
\textsf{A} & \textsf{0}\\
-\textsf{A}_{\Gamma}\textsf{R}_{\Gamma} & \textsf{A}_\Gamma
\end{array}\right]^{-1}
\left[\begin{array}{c}
\textsf{A}_\Gamma\\
\textsf{A}_\Gamma
\end{array}\right]
\end{equation}
and
\begin{equation}\label{eq:b}
\textsf{b}=(\textsf{I}_\Gamma-\textsf{R}_\Gamma
\textsf{A}^{-1} \textsf{A}_{\Gamma})\textsf{R}_{\Gamma} \textsf{A}^{-1}
\textsf F.
\end{equation}

Since the matrix \textsf{S} is not symmetric, the Schur-complement system \eqref{eq:Schur} is solved by a Krylov method, e.g., GMRES or BiCGStab \cite{Vander:2003} after computing the right-hand side \textsf{b} using formula (\ref{eq:b}) by solving separately local Stokes and Darcy problems,
as detailed in Algorithm \ref{alg:2}.

Once \textsf{b} is available, each iteration of the Krylov method requires evaluating the action of the operator \textsf{S} on a given array $\textsf g=[\textsf g_f,\, \textsf g_p]^t$ containing the DoFs of the Stokes velocity on $\Gamma_f$ and the DoFs of the Darcy pressure on $\Gamma_p$, again by solving separately local Stokes and Darcy problems; 
this procedure is outlined in Algorithms \ref{alg:3} and \ref{alg:4}. Finally, Algorithm \ref{alg:5} summarizes all the steps needed to obtain the solution of the global system \eqref{eq:algebraic-OS-icdd}.

\begin{remark1}
For the sake of clarity, Algorithm \ref{alg:1} indicates all the matrices that should be assembled to solve system (\ref{eq:Schur}). However, notice that Algorithm \ref{alg:1} can be replaced by any user's code implementing the construction of the required matrices. Moreover, in all Algorithms \ref{alg:2}--\ref{alg:5}, instead of inputting the preassembled matrices $\mathsf A_f$, $\mathsf A_{\Gamma_f}$, $\mathsf B_f$, $\mathsf C_f$, $\mathsf R_f$, $\mathsf A_p$, $\mathsf A_{\Gamma_p}$, $\mathsf B_p$, $\mathsf C_p$, and $\mathsf R_p$, one could provide instructions to execute ad-hoc software that solves the specified Stokes and Darcy problems and returns the quantities of interest $\mathsf{b}$, $\mathsf{t}$ or $\mathsf{g}$. This makes the ICDD method independent of the methodology and software used for the discretization of the Stokes and the Darcy problems. The boundary conditions imposed on the interfaces $\Gamma_f$ and $\Gamma_p$ are standard, so that they can be easily handled by any computational software.
\end{remark1}

\begin{algorithm}[h!]
\caption{Local matrices assembling (by any available code)} \label{alg:1}
\begin{algorithmic}[1]
\Procedure{Assemble\,}{meshes, data}
\State \emph{Assemble the Stokes arrays:}
\State $\textsf A_f,\, \textsf A_{\Gamma_f},\, \textsf B_f,\, \textsf C_f$ stiffness matrices
\State $\textsf f_f$ right hand side
\State $\textsf R_{\Gamma_p f}$ restriction matrix from $\overline\Omega_f\setminus\Gamma_f$ to $\Gamma_p$            
\State \emph{Assemble the Darcy arrays:}
\State $\textsf A_p,\, \textsf A_{\Gamma_p},\, \textsf B_p,\, \textsf C_p$   stiffness matrices
\State $\textsf f_p$ right hand side
\State $\textsf R_{\Gamma_f p}$ restriction matrix from $\overline\Omega_p\setminus\Gamma_p$ to $\Gamma_f$	 
\State \Return{arrays $\textsf A_f,\, \textsf A_{\Gamma_f},\, \textsf B_f,\, \textsf C_f,\,
\textsf f_f,\,  \textsf R_{\Gamma_P f},\, \textsf A_p,\, \textsf A_{\Gamma_p},\, \textsf B_p,\, \textsf C_p,\, \textsf f_p,\, \textsf R_{\Gamma_f p}$ }
\EndProcedure
\end{algorithmic}
\end{algorithm}	

\begin{algorithm}[h!]
\caption{Computation of the right-hand side $\textsf b$ of system \eqref{eq:Schur}} \label{alg:2}
\begin{algorithmic}[1]
\Procedure{SchurRHS\,}{$\textsf A_f$, $\textsf A_{\Gamma_f}$, $\textsf B_f$, $\textsf C_f$, $\textsf f_f$, $\textsf R_{\Gamma_p f}$, $\textsf A_p$,  $\textsf A_{\Gamma_p}$, $\textsf B_p$, $\textsf C_p$, $\textsf f_p$, $\textsf R_{\Gamma_f p}$}
\State\emph{Solve the primal Stokes problem with} $\widehat{\textsf{u}}_f=\textsf{0}$ \emph{on $\Gamma_f$:}
\State $ \left[\begin{array}{cc}
        \textsf A_f & \textsf B_f^t \\
        \textsf B_f & \textsf C_f
        \end{array}\right] \, 
        \left[\begin{array}{c}
        \widehat{\textsf u}_f \\
        \widehat{\textsf p}_f
        \end{array}\right] = 
        \left[\begin{array}{c}
        \textsf f_f \\
        \textsf 0
        \end{array}\right] $
\State\emph{Extract the nodal values of} $\widehat{\textsf p}_f$ \emph{on $\Gamma_p$:}
\State $\widehat{\textsf p}_{f|\Gamma_p}=\textsf R_{\Gamma_p f} \widehat{\textsf p}_f$ 
\State\emph{Solve the primal Darcy problem with} $\widehat{\textsf{p}}_p=\textsf{0}$ \emph{on $\Gamma_p$:}
\State $ \left[\begin{array}{cc}
        \textsf A_p & \textsf B_p^t \\
        \textsf B_p & \textsf C_p
        \end{array}\right] \, 
        \left[\begin{array}{c}
        \widehat{\textsf u}_p \\
        \widehat{\textsf p}_p
        \end{array}\right] = 
        \left[\begin{array}{c}
        \textsf f_p \\
        \textsf 0
        \end{array}\right] $
\State\emph{Extract the nodal values of} $\widehat{\textsf u}_p$ \emph{on $\Gamma_f$:}
\State $\widehat{\textsf u}_{p|\Gamma_f}=\textsf R_{\Gamma_f p} \widehat{\textsf u}_p$ 
\State\emph{Solve the dual Stokes problem with} $\widehat{\textsf w}_f=-\widehat{\textsf u}_{p|\Gamma_f}$ \emph{on $\Gamma_f$:}
\State $ \left[\begin{array}{cc}
        \textsf A_f & \textsf B_f^t \\
        \textsf B_f & \textsf C_f
        \end{array}\right] \, 
        \left[\begin{array}{c}
        \widehat{\textsf w}_f \\
        \widehat{\textsf q}_f
        \end{array}\right] = 
        \left[\begin{array}{c}
        \textsf A_{\Gamma_f}\widehat{\textsf u}_{p|\Gamma_f}\\
        \textsf 0 \\
        \end{array}\right] $
\State\emph{Extract the nodal values of} $\widehat{\textsf q}_f$ \emph{on $\Gamma_p$:}
\State $\widehat{\textsf q}_{f|\Gamma_p}=\textsf R_{\Gamma_p f} \widehat{\textsf q}_f$
\State\emph{Solve the dual Darcy problem with} $\widehat{\textsf q}_p=-\widehat{\textsf p}_{f|\Gamma_p}$ \emph{on $\Gamma_p$:}
\State $ \left[\begin{array}{cc}
        \textsf A_p & \textsf B_p^t \\
        \textsf B_p & \textsf C_p
        \end{array}\right] \, 
        \left[\begin{array}{c}
        \widehat{\textsf w}_p \\
        \widehat{\textsf q}_p
        \end{array}\right] = 
        \left[\begin{array}{c}
        \textsf A_{\Gamma_p}\widehat{\textsf p}_{f|\Gamma_p}\\
        \textsf 0
        \end{array}\right] $
\State\emph{Extract the nodal values of} $\widehat{\textsf w}_p$ \emph{on $\Gamma_f$:}
\State $\widehat{\textsf w}_{p|\Gamma_f}=\textsf R_{\Gamma_f p} \widehat{\textsf u}_p$ 
\State \Return{
$\textsf b=[-\widehat{\textsf u}_{p|\Gamma_f}+\widehat{\textsf w}_{p|\Gamma_f},\ \widehat{\textsf q}_{f|\Gamma_p}-\widehat{\textsf p}_{f|\Gamma_p}]^t$}
\EndProcedure
\end{algorithmic}
\end{algorithm}	

\begin{algorithm}[h!]
\caption{Given $\textsf{g}=[\textsf{g}_f,\ \textsf{g}_p]^t$, compute $\textsf t =\textsf S \textsf g$} \label{alg:3}
\begin{algorithmic}[1]
\Procedure{SchurEval\,}{$\textsf{g}=[\textsf{g}_f,\, \textsf{g}_p]^t$, $\textsf A_f$, $\textsf A_{\Gamma_f}$, $\textsf B_f$, $\textsf C_f$, $\textsf R_{\Gamma p f}$, $\textsf A_p$, $\textsf A_{\Gamma_p}$, $\textsf B_p$, $\textsf C_p$, $\textsf R_{\Gamma_f p}$} 
\State\emph{Solve the primal Stokes problem with} $\widetilde{\textsf u}_f=\textsf g_f$ \emph{on $\Gamma_f$ (all other problem data are zero):}
\State $ \left[\begin{array}{cc}
        \textsf A_f & \textsf B_f^t \\
        \textsf B_f & \textsf C_f
        \end{array}\right] \, 
        \left[\begin{array}{c}
        \widetilde{\textsf u}_f \\
        \widetilde{\textsf p}_f
        \end{array}\right] = 
        \left[\begin{array}{c}
        -\textsf A_{\Gamma_f} \textsf g_f \\
        \textsf 0
        \end{array}\right] $
\State\emph{Extract the nodal values of} $\widetilde{\textsf p}_f$ \emph{on $\Gamma_p$:}
        \State $\widetilde{\textsf p}_{f|\Gamma_p} =\textsf R_{\Gamma_p f}\widetilde{\textsf p}_f$ 
\State\emph{Solve the primal Darcy problem with $\widetilde{\mathsf p}_p=\textsf g_p$ on $\Gamma_p$ (all other problem data are zero):}
\State $ \left[\begin{array}{cc}
        \textsf A_p & \textsf B_p^t \\
        \textsf B_p & \textsf C_p
        \end{array}\right] \, 
        \left[\begin{array}{c}
        \widetilde{\textsf u}_p \\
        \widetilde{\textsf p}_p
        \end{array}\right] = 
        \left[\begin{array}{c}
        -\textsf A_{\Gamma_p} \textsf g_p \\
        \textsf 0
        \end{array}\right]$ 
\State\emph{Extract the nodal values of} $\widetilde{\mathsf u}_p$ \emph{on} $\Gamma_f$:
\State $\widetilde{\textsf u}_{p|\Gamma_f} =\textsf R_{\Gamma_f p} \widetilde{\textsf u}_p$
\State\emph{Solve the dual Stokes problem with  $\widetilde{\mathsf w}_f=\mathsf{g}_f-\widetilde{\mathsf u}_{p|\Gamma_f}$ on $\Gamma_f$:}
\State $ \left[\begin{array}{cc}
        \textsf A_f & \textsf B_f^t \\
        \textsf B_f & \textsf C_f
        \end{array}\right] \, 
        \left[\begin{array}{c}
        \widetilde{\textsf w}_f \\
        \widetilde{\textsf q}_f
        \end{array}\right] = 
        \left[\begin{array}{c}
        -\textsf A_{\Gamma_f}(\textsf{g}_f-\widetilde{\textsf u}_{p|\Gamma_f})\\
        \textsf 0 
        \end{array}\right] $
\State\emph{Extract the nodal values of} $\widetilde{\textsf q}_f$ \emph{on $\Gamma_p$:}
\State $\widetilde{\textsf q}_{f|\Gamma_p}=\textsf R_{\Gamma_p f}\, \widetilde{\textsf q}_f$
\State\emph{Solve the dual Darcy problem with  $\widetilde{\mathsf q}_p=\mathsf g_p- \widetilde{\mathsf p}_{f|\Gamma_p}$ on $\Gamma_p$}:
\State $ \left[\begin{array}{cc}
        \textsf A_p & \textsf B_p^t \\
        \textsf B_p & \textsf C_p
        \end{array}\right] \, 
        \left[\begin{array}{c}
        \widetilde{\textsf w}_p \\
        \widetilde{\textsf q}_p
        \end{array}\right] = 
        \left[\begin{array}{c}
        -\textsf A_{\Gamma p}(\textsf g_p-\widetilde{\textsf p}_{f|\Gamma_p})\\
        \textsf 0 
        \end{array}\right] $
\State\emph{Extract the nodal values of} $\widetilde{\textsf w}_p$ \emph{on $\Gamma_f$:}
\State $\widetilde{\textsf w}_{p|\Gamma_f}=\textsf R_{\Gamma_f p} \widetilde{\textsf w}_p$ 
\State \Return{
$\textsf t=[\textsf g_f-\widetilde{\textsf u}_{p|\Gamma_f}+\widetilde{\textsf w}_{p|\Gamma_f},\, \textsf g_p+\widetilde{\textsf q}_{f|\Gamma_p}-\widetilde{\textsf p}_{f|\Gamma_p}]^t$}

\EndProcedure
\end{algorithmic}
\end{algorithm}	

\begin{algorithm}[h!]
\caption{Solve the Schur complement system (\ref{eq:Schur})} \label{alg:4}
\begin{algorithmic}[1]
\Procedure{SchurSolve\,}{$\textsf{b}$, $\textsf A_f$, $\textsf A_{\Gamma_f}$, $\textsf B_f$, $\textsf C_f$, $\textsf R_{\Gamma_p f}$, $\textsf A_p$, $\textsf A_{\Gamma_p}$, $\textsf B_p$, $\textsf C_p$, $\textsf R_{\Gamma_f p}$ }
\State given $\textsf g^{(0)}=[\textsf{g}_f^{(0)},\, \textsf{g}_p^{(0)}]^t$:
\For{$k=0,\ldots,$ until convergence}
\State \emph{Krylov iteration}
\State $\ldots$
\State \emph{Evaluation of} $\textsf t^{(k)}=\textsf S \textsf g^{(k)}$:
\State $\textsf t^{(k)}$ = \Call{SchurEval\,}{$\textsf{g}^{(k)}=[\textsf{g}_f^{(k)},\, \textsf{g}_p^{(k)}]^t$, $\textsf A_f$, $\textsf A_{\Gamma_f}$, $\textsf B_f$, $\textsf C_f$, $\textsf R_{\Gamma_p f}$, $\textsf A_p$, $\textsf A_{\Gamma_p}$, $\textsf B_p$, $\textsf C_p$, $\textsf R_{\Gamma_f p}$}
\State $\ldots$
\State \emph{End Krylov iteration}
\EndFor
\State \Return{$\textsf g=[\textsf g_f, \textsf g_p]^t$}
\EndProcedure
\end{algorithmic}
\end{algorithm}

\begin{algorithm}[h!]
\caption{Solve system \eqref{eq:algebraic-OS-icdd} } \label{alg:5}
\begin{algorithmic}[1]
\Procedure{Solve}{meshes, FE, data}
\State\emph{Assemble local matrices using any available code}
\State 
$ [\textsf A_f,\, \textsf A_{\Gamma_f},\, \textsf B_f,\, \textsf C_f,\, \textsf f_f,\, \textsf R_{\Gamma_p f},\, \textsf A_p,\, \textsf A_{\Gamma_p},\, \textsf B_p,\, \textsf C_p,\, \textsf f_p,\, \textsf R_{\Gamma_f p}] \gets\Call{Assemble}{\text{meshes, FE, data}}$
\State \emph{Compute the right hand side} $\textsf b$ of the Schur--complement system
\State $\textsf b\gets\Call{SchurRHS}{\textsf A_f,\, \textsf B_f,\, \textsf C_f,\, \textsf f_f,\, \textsf R_{\Gamma_p f},\, \textsf A_p,\, \textsf B_p,\, \textsf C_p,\, \textsf f_p,\, \textsf R_{\Gamma_f p}}$
\State \emph{Solve} $\textsf S \textsf g = \textsf b$ \emph{by a Krylov method}
\State $\textsf g\gets\Call{SchurSolve}{\textsf b,\, \textsf A_f,\, \textsf A_{\Gamma_f},\, \textsf B_f,\, \textsf C_f,\, \textsf R_{\Gamma_p f},\, \textsf A_p,\, \textsf A_{\Gamma_p},\, \textsf B_p,\, \textsf C_p,\, \textsf R_{\Gamma_f p}}$
\State\emph{Solve the complete Stokes problem:}
\State 
       $ \left[\begin{array}{cc}
        \textsf A_f & \textsf B_f^t \\
        \textsf B_f & \textsf C_f
        \end{array}\right] \, 
        \left[\begin{array}{c}
        {\textsf u}_f \\
        {\textsf p}_f
        \end{array}\right] = 
        \left[\begin{array}{c}
        \textsf f_f -\textsf A_{\Gamma_f} \textsf g_f \\
        \textsf 0
        \end{array}\right] $
\State\emph{Solve the complete Darcy problem:}
\State 
       $ \left[\begin{array}{cc}
        \textsf A_p & \textsf B_p^t\\
        \textsf B_p & \textsf C_p
        \end{array}\right] \, 
        \left[\begin{array}{c}
        {\textsf u}_p \\
        {\textsf p}_p
        \end{array}\right] = 
        \left[\begin{array}{c}
        \textsf f_p -\textsf A_{\Gamma_p} \textsf g_p \\
        \textsf 0
        \end{array}\right]$ 
\State \Return{$\textsf u_f,\, \textsf p_f,\, \textsf u_p,\, \textsf p_p$}
\EndProcedure
\end{algorithmic}
\end{algorithm}	

All the numerical tests involving the ICDD method reported in this work have been performed using the BiCGStab method, and at most 5 iterations were needed to reduce the norm of the initial residual of the Schur-complement system \eqref{eq:Schur} by 8 orders of magnitude. For a detailed study of the convergence properties of the ICDD method for Stokes-Darcy, we refer to \cite{Discacciati:2016:SINUM}.

\smallskip

Finally, we observe that problem \eqref{eq:stokesFluid}--\eqref{eq:continuityVelocity} could be solved by the Alternating Schwarz Method (ASM) \cite{Quarteroni:1999, Toselli:2005}. However, we warn the reader that the convergence rate of ASM strongly depends on the overlap thickness, say, $\delta$, with the number of iterations needed to reach convergence typically behaving like $\delta^{-2}$. On the contrary, solving the Schur-complement system \eqref{eq:Schur} is a very robust approach and the convergence rate is independent of $\delta$ provided that the permeability $\textbf K$ is small enough, as it is the case in real physical situations.

\section{Physical interpretation of the Stokes-Darcy model with overlap}
\label{sec:microscale}

In this section, we provide a physical interpretation of the Stokes-Darcy model with overlap described in Sect.~\ref{sec:stokesDarcy} and we characterize the interfaces $\Gamma_f$ and $\Gamma_p$ where the matching conditions \eqref{eq:continuityVelocity} and \eqref{eq:continuityPressure} are imposed, respectively. To this aim, we introduce a reference microscale Stokes problem that we will also use for numerical validation in Sect.~\ref{sec:validation}.

\subsection{Reference microscale model}\label{sec:microscaleModel}

\begin{figure}[bht]
\begin{center}
 	\begin{tikzpicture}[scale=0.75]
       \fill[gray!12] (-0.5,2.2) rectangle (4.5,2.8);
        \fill[gray!40] (-0.25, -0.25) circle (4pt);
        \fill[gray!40] ( 0.25, -0.25) circle (4pt);
        \fill[gray!40] ( 0.75, -0.25) circle (4pt);
        \fill[gray!40] ( 1.25, -0.25) circle (4pt);
        \fill[gray!40] ( 1.75, -0.25) circle (4pt);
        \fill[gray!40] ( 2.25, -0.25) circle (4pt);
        \fill[gray!40] ( 2.75, -0.25) circle (4pt);
        \fill[gray!40] ( 3.25, -0.25) circle (4pt);
        \fill[gray!40] ( 3.75, -0.25) circle (4pt);
        \fill[gray!40] ( 4.25, -0.25) circle (4pt);
        \fill[gray!40] (-0.25, 0.25) circle (4pt);
        \fill[gray!40] ( 0.25, 0.25) circle (4pt);
        \fill[gray!40] ( 0.75, 0.25) circle (4pt);
        \fill[gray!40] ( 1.25, 0.25) circle (4pt);
        \fill[gray!40] ( 1.75, 0.25) circle (4pt);
        \fill[gray!40] ( 2.25, 0.25) circle (4pt);
        \fill[gray!40] ( 2.75, 0.25) circle (4pt);
        \fill[gray!40] ( 3.25, 0.25) circle (4pt);
        \fill[gray!40] ( 3.75, 0.25) circle (4pt);
        \fill[gray!40] ( 4.25, 0.25) circle (4pt);
        \fill[gray!40] (-0.25, 0.75) circle (4pt);
        \fill[gray!40] ( 0.25, 0.75) circle (4pt);
        \fill[gray!40] ( 0.75, 0.75) circle (4pt);
        \fill[gray!40] ( 1.25, 0.75) circle (4pt);
        \fill[gray!40] ( 1.75, 0.75) circle (4pt);
        \fill[gray!40] ( 2.25, 0.75) circle (4pt);
        \fill[gray!40] ( 2.75, 0.75) circle (4pt);
        \fill[gray!40] ( 3.25, 0.75) circle (4pt);
        \fill[gray!40] ( 3.75, 0.75) circle (4pt);
        \fill[gray!40] ( 4.25, 0.75) circle (4pt);
        \fill[gray!40] (-0.25, 1.25) circle (4pt);
        \fill[gray!40] ( 0.25, 1.25) circle (4pt);
        \fill[gray!40] ( 0.75, 1.25) circle (4pt);
        \fill[gray!40] ( 1.25, 1.25) circle (4pt);
        \fill[gray!40] ( 1.75, 1.25) circle (4pt);
        \fill[gray!40] ( 2.25, 1.25) circle (4pt);
        \fill[gray!40] ( 2.75, 1.25) circle (4pt);
        \fill[gray!40] ( 3.25, 1.25) circle (4pt);
        \fill[gray!40] ( 3.75, 1.25) circle (4pt);
        \fill[gray!40] ( 4.25, 1.25) circle (4pt);
        \fill[gray!40] (-0.25, 1.75) circle (4pt);
        \fill[gray!40] ( 0.25, 1.75) circle (4pt);
        \fill[gray!40] ( 0.75, 1.75) circle (4pt);
        \fill[gray!40] ( 1.25, 1.75) circle (4pt);
        \fill[gray!40] ( 1.75, 1.75) circle (4pt);
        \fill[gray!40] ( 2.25, 1.75) circle (4pt);
        \fill[gray!40] ( 2.75, 1.75) circle (4pt);
        \fill[gray!40] ( 3.25, 1.75) circle (4pt);
        \fill[gray!40] ( 3.75, 1.75) circle (4pt);
        \fill[gray!40] ( 4.25, 1.75) circle (4pt);
        \fill[gray!40] (-0.25, 2.25) circle (4pt);
        \fill[gray!40] ( 0.25, 2.25) circle (4pt);
        \fill[gray!40] ( 0.75, 2.25) circle (4pt);
        \fill[gray!40] ( 1.25, 2.25) circle (4pt);
        \fill[gray!40] ( 1.75, 2.25) circle (4pt);
        \fill[gray!40] ( 2.25, 2.25) circle (4pt);
        \fill[gray!40] ( 2.75, 2.25) circle (4pt);
        \fill[gray!40] ( 3.25, 2.25) circle (4pt);
        \fill[gray!40] ( 3.75, 2.25) circle (4pt);
        \fill[gray!40] ( 4.25, 2.25) circle (4pt);
        %
        \draw[dotted] (-0.5,0)--(4.5,0);
        \draw[dotted] (-0.5,0.5)--(4.5,0.5);
        \draw[dotted] (-0.5,1)--(4.5,1);
        \draw[dotted] (-0.5,1.5)--(4.5,1.5);
        \draw[dotted] (-0.5,2)--(4.5,2);
        \draw[dotted] (-0.5,2.5)--(4.5,2.5);
        \draw[dotted] (0,-0.5)--(0,2.5);
        \draw[dotted] (0.5,-0.5)--(0.5,2.5);
        \draw[dotted] (1,-0.5)--(1,2.5);
        \draw[dotted] (1.5,-0.5)--(1.5,2.5);
        \draw[dotted] (2,-0.5)--(2,2.5);
        \draw[dotted] (2.5,-0.5)--(2.5,2.5);
        \draw[dotted] (3,-0.5)--(3,2.5);
        \draw[dotted] (3.5,-0.5)--(3.5,2.5);
        \draw[dotted] (4,-0.5)--(4,2.5);
        %
        \draw[dashed] (4.0,1.5)--(8.0,2.2);
        \node [black] at (7.2,1.5) {{\footnotesize{$\times$}}$\ell$};
        \draw[<-] (6.5,1.25)--(7.7,1.25);
        \draw[dashed] (4.0,1.0)--(8.0,0.3);
        \draw (8.0,0.3) rectangle (9.9,2.2);
        \fill[gray!40] (8.95,1.25) circle (13pt);
        \node [black] at (9.0,2.5) {$Y$};
        \node [black] at (8.95,1.25) {$Y_s$};
        \node [black] at (9.5,0.6) {$Y_f$};
        \fill [black] (8.0,0.3) circle (2pt);
        \node [black] at (8.0,-0.1) {$(0,0)$};
        \fill [black] (9.9,2.2) circle (2pt);
        \node [black] at (10.5,2.2) {$(1,1)$};
        %
        \draw (-0.5,-0.5)--(4.5,-0.5);
        \draw (-0.5,-0.5)--(-0.5,5.5);
        \draw (-0.5,5.5)--(4.5,5.5);
        \draw (4.5,-0.5)--(4.5,5.5);
        \draw [dashed] (-0.5,2.2)--(4.5,2.2);
        \draw [dashed] (-0.5,2.8)--(4.5,2.8);
        \node [black] at (5.2,2.5) {$(\frac{L}{2},0)$};
        \fill [black] (4.5, 2.5) circle (2pt);
        
        \fill[gray!40] ( 4.25, 2.25) circle (4pt);

 		\node[black] at (-1.1, 1) {$\Omega_{\varepsilon p}$};
 		\draw[<->, black,dotted]  (-0.7, -0.5) -- (-0.7, 2.49);
		
 		\node[black] at (-1.1, 4) {$\Omega_{\varepsilon f}$};
 		\draw[<->, black,dotted]  (-0.7, 2.56) -- (-0.7, 5.5);
				
 		\node[black] at (2.2, 3.1) {$\Omega_{\varepsilon t} \approx O(\ell)$};

        \node[black] at (2, 4.3) {$L$};
        \draw[<->, black,dotted]  (-0.5, 4.05) -- (4.5, 4.05);
        
        \node[black] at (-0.8, 6.0) {$\Omega_\varepsilon$};
        
        \node[black] at (2, 5.9) {$\Gamma_{top}$};
        \draw[<->, black,dotted]  (-0.5, 5.6) -- (4.5, 5.6);
        \fill [black] (4.5, 5.5) circle (2pt);
        \node [black] at (5.3,5.5) {$(\frac{L}{2},H)$};
 		
		\node[black] at (4.9, 1.25) {$Y^\ell$};
        \node[black] at (5.0, 0.25) {$\ell$};
        \draw[<->, black,dotted]  (4.7, 0.0) -- (4.7, 0.5);
        
        \node[black] at (2.2, -0.9) {$\Gamma_{bottom}$};
        \draw[<->, black,dotted]  (-0.5, -0.65) -- (4.5, -0.65);
        
        \fill [black] (-0.5, -0.5) circle (2pt);
        \node [black] at (-0.6,-1.0) {$(-\frac{L}{2},-d)$};
 	\end{tikzpicture}
	\end{center}
  \caption{Schematic representation of the microscale domain $\Omega_\varepsilon$.}
  \label{fig:domainMicroscopic}
\end{figure}
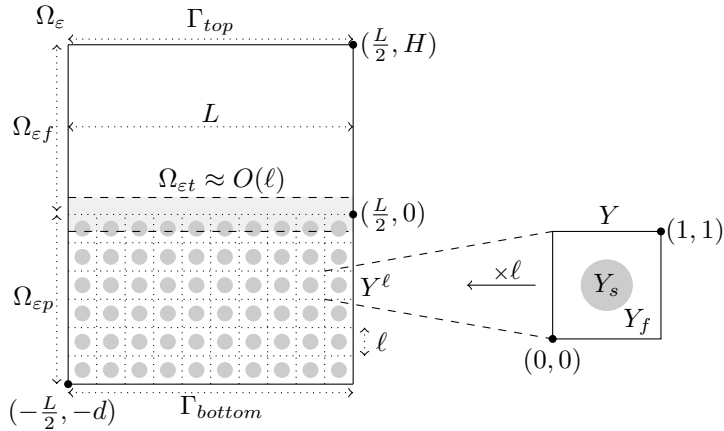

Let $\Omega_\varepsilon$ be a bounded domain filled by an incompressible viscous fluid and partially formed by a periodic array of solid obstacles to represent an isotropic porous medium, see Fig.~\ref{fig:domainMicroscopic}. More precisely, let $\Omega_\varepsilon = \Omega_{\varepsilon f} \cup \Omega_{\varepsilon p}$ with $\Omega_{\varepsilon f} = (-\frac{L}{2},\frac{L}{2}) \times (0,H)$ and $\Omega_{\varepsilon p} = (-\frac{L}{2},\frac{L}{2}) \times (-d,0)$ with $L,H,d>0$, and let $\Omega_{\varepsilon p}$ contain a periodic arrangement of solid obstacles. In particular, let $Y=(0,1)^2$ be a dimensionless unit cell including a solid obstacle $Y_s$ centred at $(0.5,0.5)$, while let $Y_f = Y\setminus Y_s$ be the fluid region around the obstacle. The cell $Y$ is scaled by a suitable characteristic microscale length $0<\ell<1$, and the resulting scaled cell $Y^\ell$ is periodically repeated to cover the whole region $\Omega_{\varepsilon p}$ as shown in Fig.~\ref{fig:domainMicroscopic}. The shape of the solid obstacles $Y_s$ is chosen to guarantee that the porous medium represented by the domain $\Omega_{\varepsilon p}$ is isotropic, i.e., the symmetric and positive definite permeability tensor $\mathbf{K}$ that characterizes the porous material is $\mathbf{K}=K\mathbf{I}$ with $K>0$ constant. We assume that separation of scales holds between the free-fluid region $\Omega_{\varepsilon f}$ and the region $\Omega_{\varepsilon p}$ occupied by the solid inclusions, i.e.
\begin{equation*}
 \varepsilon = \frac{\ell}{L} \ll 1\,,
\end{equation*}
where $L$ is the characteristic macroscale length in $\Omega_{\varepsilon f}$.

Under the hypothesis that the Reynolds number in $\Omega_\varepsilon$ is small, the fluid can be modelled by the dimensional Stokes equations
\begin{subequations}\label{eq:stokesGlobal}
  \begin{eqnarray}
    -\mu \, \Delta \mathbf{u}_\varepsilon + \nabla p_\varepsilon = \mathbf{f} && \text{in } \Omega_\varepsilon \, ,\\
                                      \nabla \cdot \mathbf{u}_\varepsilon = 0 && \text{in } \Omega_\varepsilon \, ,
  \end{eqnarray}
\end{subequations}
where $\mu$ and $\textbf{f}$ are as in \eqref{eq:stokesFluid_momentum},
$\mathbf{u}_\varepsilon$ is the fluid velocity 
and $p_\varepsilon$ is the pressure.           
(Suitable boundary conditions to ensure the well-posedness of \eqref{eq:stokesGlobal} will be specified when needed.)

\begin{figure}[bht]
\begin{center}
\begin{tabular}{lcr}
\includegraphics[valign=t,width=0.25\textwidth]{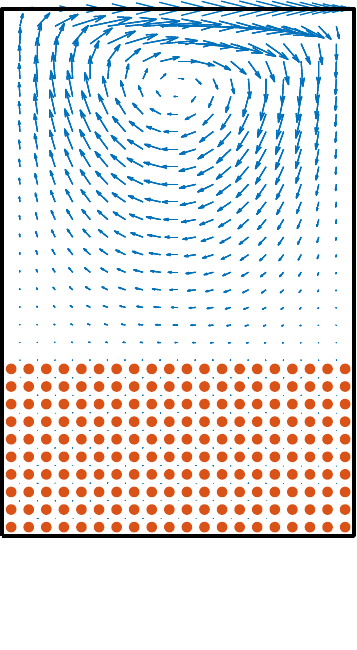}\hskip 1.cm
\includegraphics[valign=t,width=0.13\textwidth]{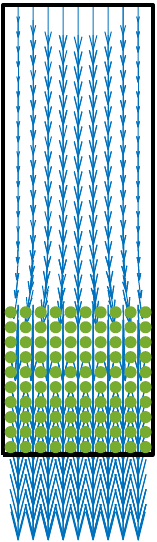}\hskip 1.cm
\includegraphics[valign=t,width=0.25\textwidth]{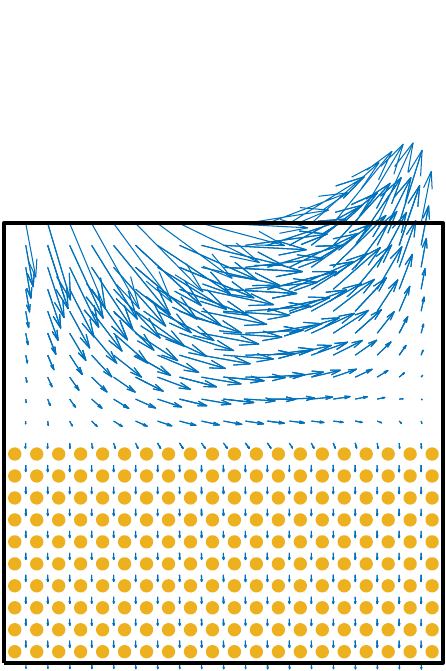}
\end{tabular}
\end{center}
 \caption{Computational domain and velocity field for the configurations with pressure gradient (almost) parallel (Test case \#1, left), perpendicular (Test case \#2, centre), and oblique (Test case \#3, right) to the porous medium.}
  \label{fig:test1-2-3}
\end{figure}

The assumption of low Reynolds number implies that pressure forces are balanced by viscous forces \cite[Sect.~9.6]{Kundu:2012} in $\Omega_\varepsilon$, and we assume that the pressure $p_\varepsilon$ is bounded in the neighbourhood of $\Omega_{\varepsilon f} \cap \Omega_{\varepsilon p}$.
While the microscopic velocity $\mathbf{u}_\varepsilon$ and pressure $p_\varepsilon$ remain continuous in $\Omega_\varepsilon$, their magnitude may significantly vary across the domain, and it is well understood (see, e.g., \cite{Ene:1975:JM,Levy:1975:IJES,LeBars:2006:JFM,Lacis:2017:JFM,Bottaro:2019:JFM,Naqvi:2021:IJMF}) that there exists an intermediate transition region $\Omega_{\varepsilon t}$ of characteristic thickness $\ell$  (see Fig.~\ref{fig:domainMicroscopic}) between the free-fluid region $\Omega_{\varepsilon f}$ and the bulk region of the porous medium $\Omega_{\varepsilon p}$ where the fluid rapidly changes its behaviour. Moreover, the characteristic magnitude of the velocity and pressure varies significantly depending on the orientation of the flow near the porous surface as already observed in \cite{Levy:1975:IJES}. To study this behaviour at the microscale, we consider three configurations where the pressure gradient is either arbitrary (namely, parallel or oblique) or perpendicular to the porous medium as illustrated in Fig.~\ref{fig:test1-2-3}.

The case when the gradient of the pressure field $p_\varepsilon$ is parallel to the porous medium has been extensively studied in the literature, see, e.g., \cite{Beavers:1967:JFM,Carraro:2013:JFM,Jager:1996:AnnPisa,Jager:2000:SIAM,Jager:2001:SISC,Lacis:2017:JFM,LeBars:2006:JFM,Levy:1975:IJES,Naqvi:2021:IJMF}. 

In general, if the pressure gradient is either parallel or oblique but not perpendicular to the porous medium, in the fluid region $\Omega_{\varepsilon f}$ the balance of viscous and pressure forces implies that the characteristic pressure is $P_f=\mu\,U_f/L$, with $U_f$ and $L$ being the characteristic velocity and length in $\Omega_{\varepsilon f}$. In the transition region $\Omega_{\varepsilon t}$, the velocity decreases due to the presence of solid obstacles and it is generally agreed that therein the characteristic velocity is
\begin{equation}\label{eq:estimateVelTransition}
  U_t = \varepsilon\,U_f
\end{equation}
(see, e.g., \cite{Bottaro:2020:Meccanica,Naqvi:2021:IJMF}). Therefore, taking $\ell$ as the characteristic length in $\Omega_{\varepsilon t}$, the balance of viscous and pressure forces (i.e., $P_t=\mu\,U_t/\ell$) implies that in $\Omega_{\varepsilon t}$ the characteristic pressure is
\begin{equation}\label{eq:pressureTransition}
P_t=P_f\,,
\end{equation}
i.e., the pressure remains of the same order of magnitude as in $\Omega_{\varepsilon f}$. Inside the porous medium domain $\Omega_{\varepsilon p}$, 
the pressure gradients at the macroscale 
are of the same order as those at the macroscale in $\Omega_{\varepsilon f}$ close to the porous medium, 
i.e.,
$\frac{P_f}{L} = \frac{P_p}{L}$, 
where $P_p$ is the macroscale characteristic pressure in $\Omega_{\varepsilon p}$ \cite[Sect.~7.1]{Ene:1975:JM}. As a consequence of the balance of forces and in view of the fact that the pressure in the porous domain is governed by the macroscale, it holds $\mu\frac{U_p}{\ell^2}=\frac{P_p}{L}=\mu\frac{U_f}{L^2}$, so that the characteristic fluid velocity in $\Omega_{\varepsilon p}$ is
\begin{equation}\label{eq:estimateVelPorous}
  U_p = \varepsilon^2\,U_f\,.
\end{equation}
Therefore, as pointed out in the seminal work \cite{Levy:1975:IJES}, in the case of non-perpendicular pressure gradient, there holds $U_f \gg U_p$, while the pressure approximately remains of the same order of magnitude inside the whole domain so that the continuity of pressure can be considered as a first-order matching condition to relate the fluid regimes in $\Omega_{\varepsilon f}$ and in $\Omega_{\varepsilon p}$.

\smallskip

On the other hand, if the gradient of the pressure $p_\varepsilon$ is perpendicular to the porous medium, \cite{Levy:1975:IJES} observed that the characteristic velocities $U_f$ and $U_p$ are both of the same order of magnitude $\sim\varepsilon^2$ in $\Omega_\varepsilon$. The normal component of the velocity is larger than the tangential one in the free-fluid domain, they are of the same order of magnitude inside the transition layer, then the difference is again evident in the bulk region of the porous medium where the tangential component rapidly decreases.
Moreover, the pressure undergoes a much larger gradient inside $\Omega_{\varepsilon p}$ than in $\Omega_{\varepsilon f}$ but, to the first order, it remains continuous in $\Omega_\varepsilon$ and almost constant along the upper row of solid inclusions in $\Omega_\varepsilon$. The continuity of the pressure was also noticed in \cite{Carraro:2013:JFM,Carraro:2015:CMAME} in the presence of an isotropic porous medium.

\smallskip

\begin{figure}[bht]
  \centering
  \includegraphics[width=0.28\textwidth]{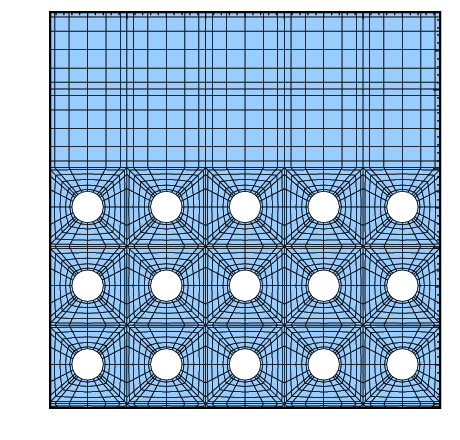}
  \caption{Example of computational mesh for the spectral element discretization of \eqref{eq:stokesGlobal}.
  }
  \label{fig:mesh}
\end{figure}

These observations are also confirmed by the numerical results at the microscale that we performed considering the three configurations of Fig.~\ref{fig:test1-2-3} described below. In all three test cases, the fluid has dynamic viscosity $\mu=10^{-3}$ kg/(m\,s) and density $\rho=10^3$ kg/m$^3$. The Stokes problem \eqref{eq:stokesGlobal} is solved using stabilized SEM \cite{Gervasio:1998:NMPDE} as in Sect.~\ref{sec:ICDD} with local polynomial degree 6 in $\Omega_\varepsilon$ on a structured mesh of size $h \approx \ell$ in $\Omega_{\varepsilon p}$ and $h = \max (\ell,\frac{1}{20})$ in $\Omega_{\varepsilon f}$. An example of the computational mesh on a subregion of $\Omega_\varepsilon$ is shown in Fig.~\ref{fig:mesh}.

\medskip

\paragraph{Test case \#1: Lid--driven cavity}\mbox{}\label{sec:test1}

We consider a lid--driven cavity setting similar to the one proposed in \cite{Lacis:2017:JFM} for which the flow is almost parallel to the upper row of solid obstacles of $\Omega_{\varepsilon p}$. In particular, the computational domain is $\Omega_\varepsilon = (-0.5,0.5)\times (-0.5,1)$~m ($L=1$~m, see also Fig.~\ref{fig:domainMicroscopic} for the notation) and we impose a null external force $\textbf{f}=\mathbf{0}$ and boundary conditions  $\mathbf{u}_\varepsilon = ((1-4x^2)10^{-6},0)$~m/s at $\Gamma_{top} = (-0.5,0.5) \times \{1 \}$~m and $\mathbf{u}_\varepsilon = \mathbf{0}$ on $\partial\Omega_\varepsilon \setminus \Gamma_{top}$, together with the constraint $\int_{\Omega_\varepsilon} p_\varepsilon = 0$ to ensure the well-posedness of the Stokes problem \eqref{eq:stokesGlobal}. The Dirichlet boundary condition on $\Gamma_{top}$ makes the characteristic velocity be $U_f=10^{-6}$~m/s so that the Reynolds number is $Re=L \,\rho \,U_f/\mu=1.$ 

\begin{figure}[bht]
  \hspace*{-6mm}
  \begin{tabular}{ccc}
  $u_{\varepsilon,1}/U_f$ & $u_{\varepsilon,2}/U_f$ & $p_\varepsilon/P_f$ \\
  \includegraphics[width=0.33\textwidth]{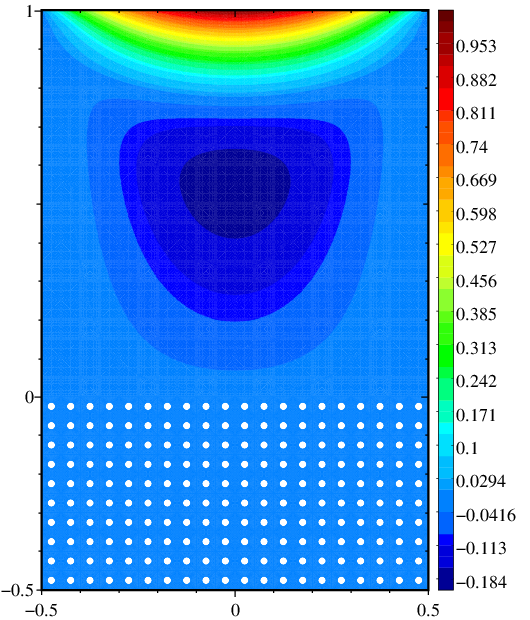} &
  \includegraphics[width=0.33\textwidth]{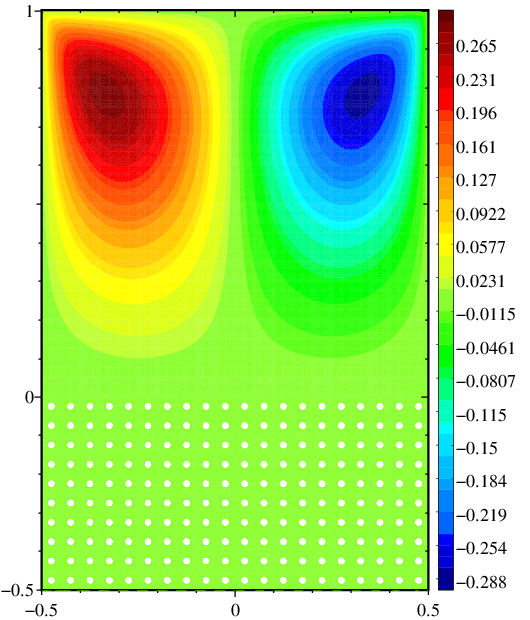} &
  \includegraphics[width=0.33\textwidth]{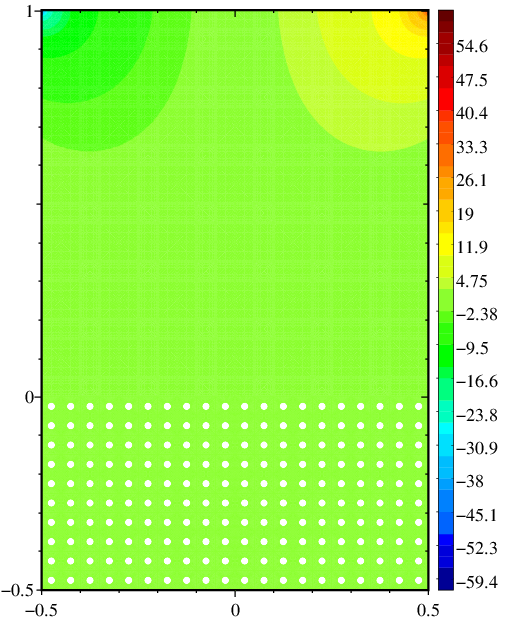}\\
  \includegraphics[width=0.33\textwidth]{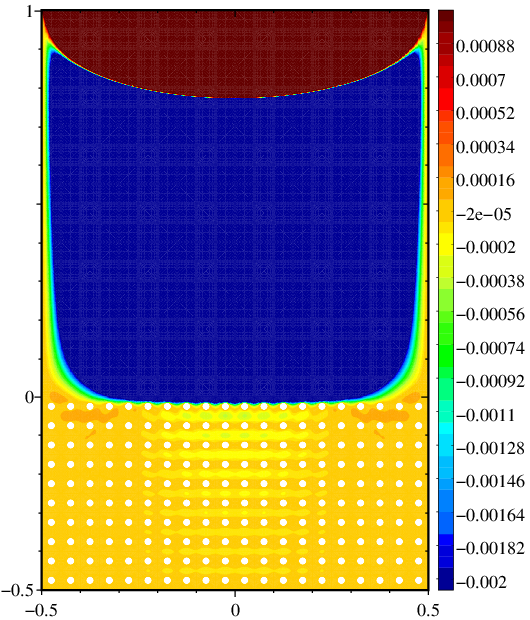} &
  \includegraphics[width=0.33\textwidth]{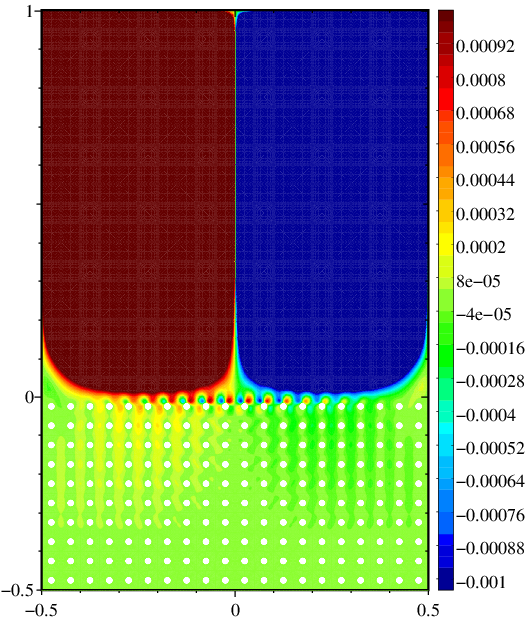} &
  \includegraphics[width=0.33\textwidth]{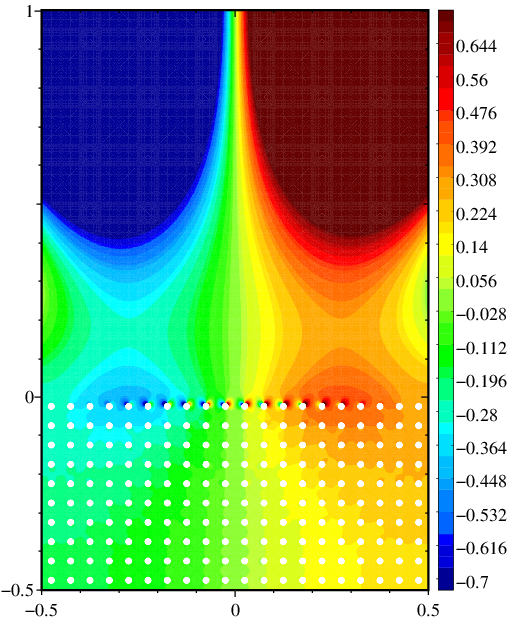}
  \end{tabular}
  \caption{\emph{Test case \#1}. From left to right: horizontal and vertical components of the dimensionless microscale velocity and pressure with $\ell=\frac{1}{20}$ and circular obstacles with $r=0.2\,\ell$. (Top) The full solution where only the macroscale behaviour is visible; (bottom) solutions plotted with a different colour scale
  to highlight small variations both in the transition region and in the porous region.}
  \label{fig:cavity}
\end{figure}

\begin{figure}[bht]
  \centering
  \includegraphics[width=\textwidth]{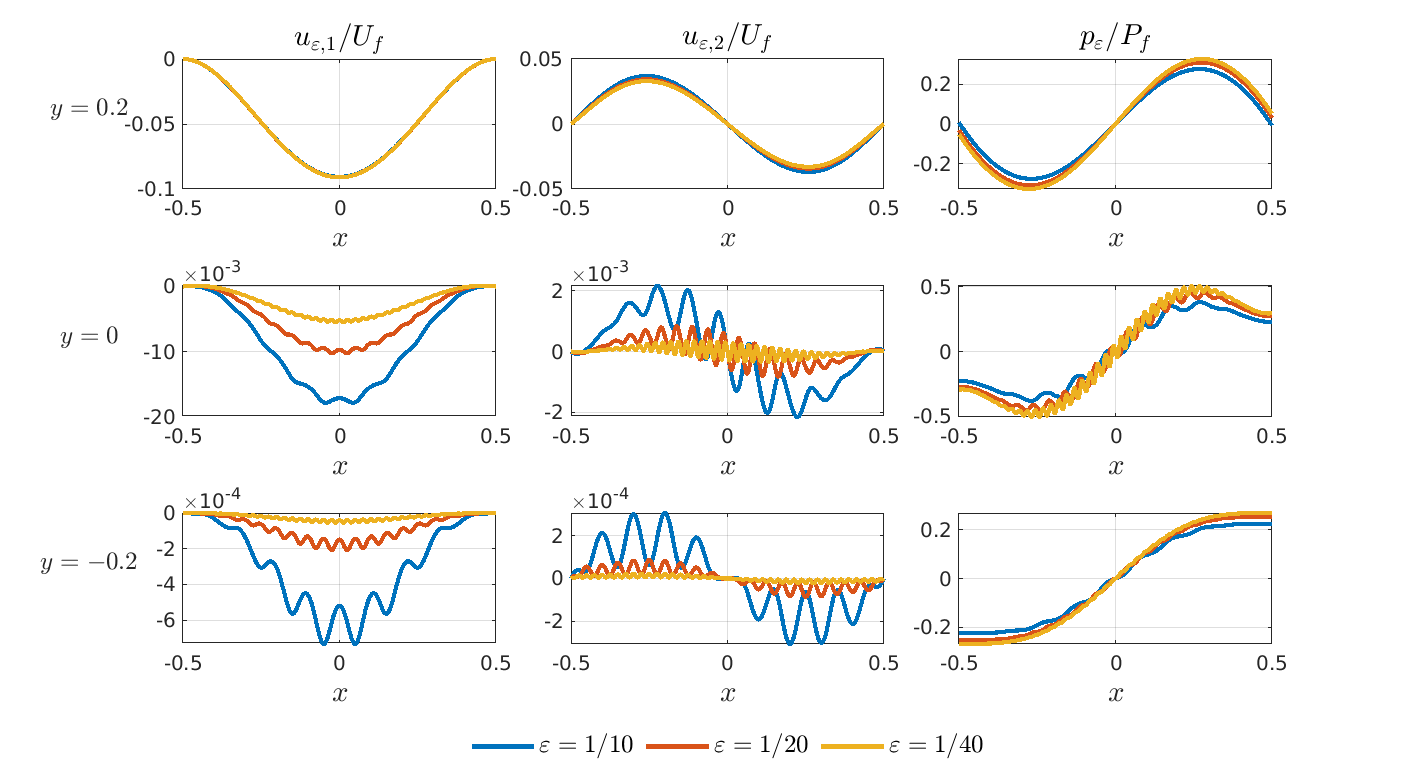}
  \caption{\emph{Test case \#1}. Velocity and pressure profiles of the dimensionless microscale solution at three vertical locations and for three different values of $\varepsilon=\ell$. Circular obstacles with $r=0.2\,\ell$.}
  \label{fig:cavitySections}
\end{figure}

The velocity and pressure computed in $\Omega_\varepsilon$ with circular obstacles characterized by $\ell=\frac{1}{20}$ and radius $r=0.2\,\ell$ are shown in Fig.~\ref{fig:cavity}, while Fig.~\ref{fig:cavitySections} reports the dimensionless velocity and pressure profiles for $\ell = \frac{1}{10}, \, \frac{1}{20}, \, \frac{1}{40}$ and $r=0.2\,\ell$ at different vertical locations: $y=0.2$ (in the free-fluid region), $y=0$ (on top of the first row of solid inclusions) and $y=-0.2$ (among the solid obstacles).

The numerical results for Test case \#1 confirm that both components of the velocity are continuous in $\Omega_\varepsilon$ and that they undergo a significant change of regime especially in the neighbourhood of the top row of solid inclusions. This is visible both in Fig.~\ref{fig:cavity} (bottom row) and in Fig.~\ref{fig:cavitySections} where it can be seen that the two components of the velocity decrease by up to five orders of magnitude moving from $\Omega_{\varepsilon f}$ to $\Omega_{\varepsilon p}$. Figure~\ref{fig:cavitySections} also shows that when $\ell$ is halved, the magnitude of the velocity at $y=0.2$ (in $\Omega_{\varepsilon f}$) does not change, while this is approximately divided by $2$ at $y=0$ (in $\Omega_{\varepsilon t}$) and by $4$ at $y=-0.2$ (in $\Omega_{\varepsilon p}$). Notice that in this test case $L=1$ m and $\ell=\varepsilon$, so that the numerical results validate estimates \eqref{eq:estimateVelTransition} and \eqref{eq:estimateVelPorous}, see also, e.g., \cite{Beavers:1967:JFM,Ene:1975:JM,Levy:1975:IJES,Neale:1974:CJCE}. In particular, since only at the Darcy regime the velocity scales like $\varepsilon^2$ (see, e.g., \cite{Allaire:1989:AA,Mei:2010}), similarly to \cite{LeBars:2006:JFM,Neale:1974:CJCE} we can also conclude that the Stokes regime remains valid until a depth proportional to $\ell$ below the ideal upper surface that delimits the porous region (in our case, $y=0$).

As concerns the pressure, significant fluctuations can be observed around the solid obstacles of the upper row at length scale $\ell$, as also pointed out, e.g., in \cite{Carraro:2013:JFM}, and they become more localized around the obstacles as $\varepsilon \to 0$. However, the pressure remains overall continuous in the domain, in agreement with theoretical results (see, e.g., \cite{LeBars:2006:JFM,Levy:1975:IJES,Neale:1974:CJCE}) and with numerical experiments performed for isotropic porous media \cite{Carraro:2013:JFM}. Moreover, the pressure maintains the same order of magnitude in $\Omega_{\varepsilon f}$ and in the transition region $\Omega_{\varepsilon t}$ (see also \eqref{eq:pressureTransition}).

\medskip

\paragraph{Test case \#2: Normal forced filtration}\mbox{}

We consider the computational domain $\Omega_\varepsilon = (-0.25,0.25)\times (-0.5,1)$~m ($L=0.5$~m) 
filled by solid obstacles in the lower region 
$\Omega_{\varepsilon p} = (-0.25,0.25) \times (-0.5,0)$~m.
We impose  null external force $\mathbf{f}=\mathbf{0}$, homogeneous Dirichlet boundary conditions $\mathbf{u}_\varepsilon = \mathbf{0}$~m/s on the vertical lateral sides, homogeneous normal stress $(\mu\nabla\mathbf{u}_\varepsilon - p_\varepsilon \mathbf{I})\,\mathbf{n}=\mathbf{0}$~kg/(m\,s$^2$) at 
$\Gamma_{top} = (-0.25,0.25) \times \{1 \}$~m,
and $(\mu\nabla\mathbf{u}_\varepsilon - p_\varepsilon\mathbf{I})\,\mathbf{n}=(0,-10^{-7})$~kg/(m\,s$^2$) at $\Gamma_{bottom}=(-0.25,0.25)\times\{-0.5\}$~m. 
This set of data makes the characteristic velocity be $U_f\sim 10^{-6}$~m/s and $Re\sim 1$.
The fluid is forced to move from the top to the bottom of the domain and the horizontal component of the velocity is almost null in the whole domain, see Fig.~\ref{fig:test1-2-3} (middle).

\begin{figure}[h!]
  \hspace*{-6mm}
  \begin{tabular}{ccc}
   $u_{\varepsilon,1}/U_f$ & $u_{\varepsilon,2}/U_f$ & $p_\varepsilon/P_f$ \\
  \includegraphics[width=0.33\textwidth]{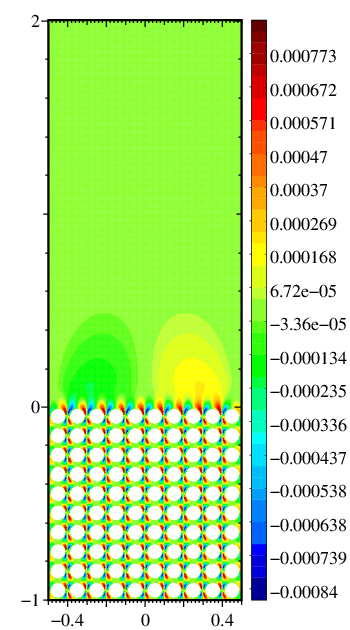} &
  \includegraphics[width=0.33\textwidth]{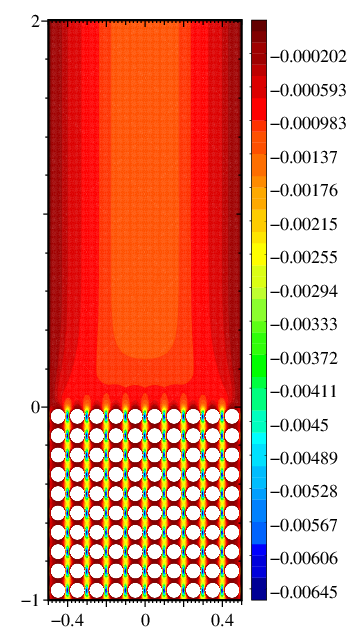} &
  \includegraphics[width=0.33\textwidth]{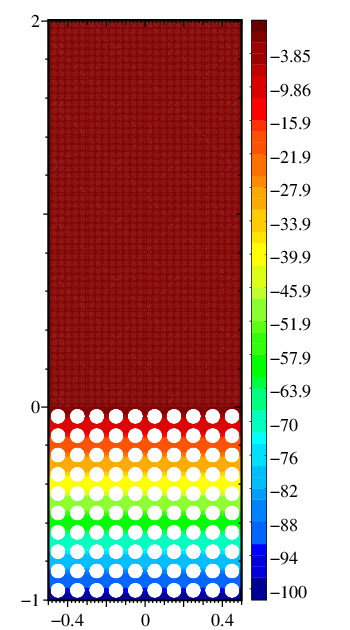}
  \end{tabular}
  \caption{\emph{Test case \#2}. From left to right: horizontal and vertical components of the dimensionless microscale velocity and pressure with $\ell=\frac{1}{20}$ and circular obstacles with $r=0.4\,\ell$.}
  \label{fig:nff_stokes_circle}
  \end{figure}

\begin{figure}[bht]
 \centering
 \includegraphics[width=\textwidth]{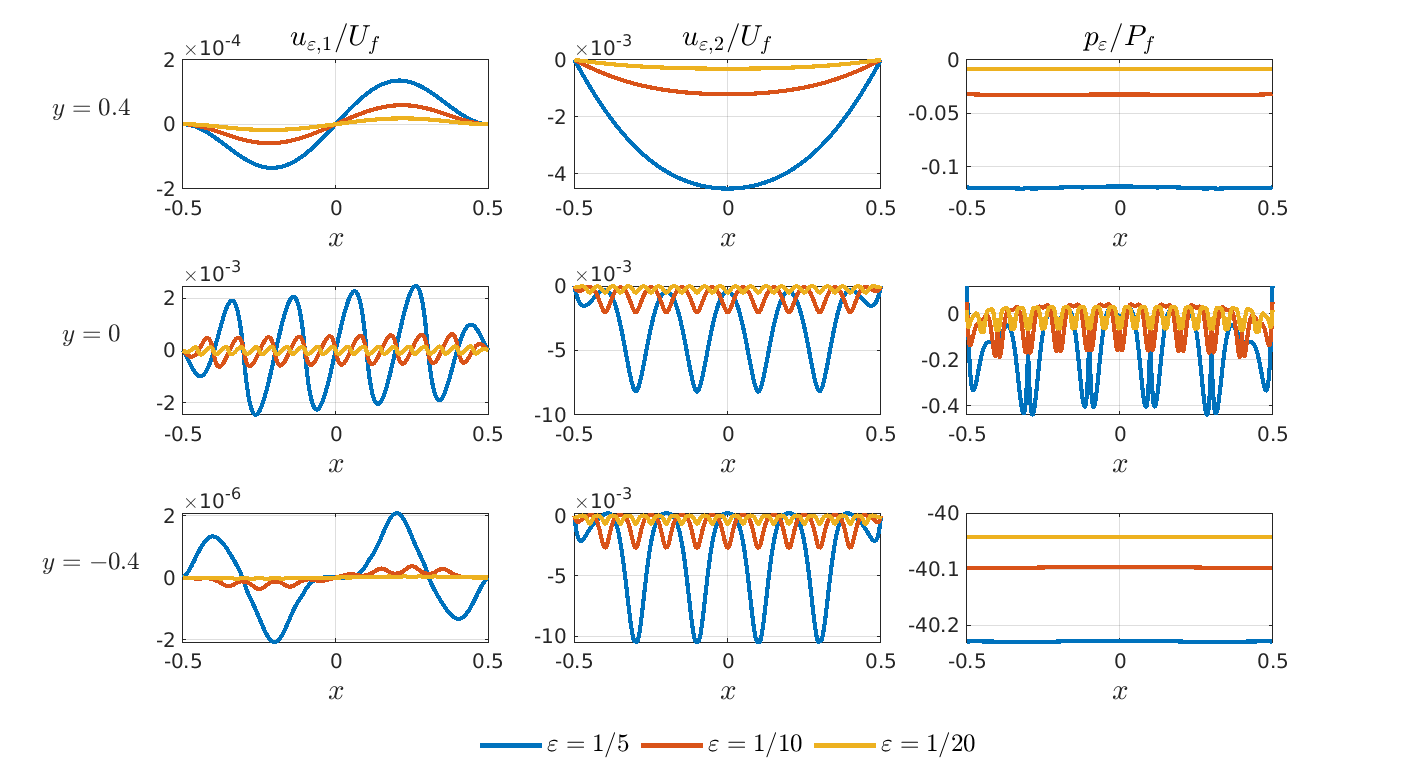}
 \caption{\emph{Test case \#2}. Velocity and pressure profiles of the dimensionless microscale solution at three vertical locations and for three different values of $\varepsilon=2\,\ell$. Circular obstacles with $r=0.4\,\ell$.}
 \label{fig:nff_sections_epsilon}
\end{figure}

The velocity and pressure computed in $\Omega_\varepsilon$ with circular obstacles characterized by $\ell=\frac{1}{20}$ and radius $r=0.4\,\ell$ are shown in Fig.~\ref{fig:nff_stokes_circle}, while Fig.~\ref{fig:nff_sections_epsilon} reports the dimensionless velocity and pressure profiles for $\ell = \frac{1}{10}, \, \frac{1}{20}, \, \frac{1}{40}$ (equivalently, $\varepsilon = \frac{1}{5}, \, \frac{1}{10}, \, \frac{1}{20}$) and $r=0.4\,\ell$ at three different vertical locations: $y=0.4$ (in the free-fluid region), $y=0$ (on top of the first row of solid inclusions), and $y=-0.4$ (among the solid obstacles). (Remark that in the dimensional setting, the three vertical coordinates coincide with those of Test case \#1 and Test case \#3 since here $L=0.5$~m instead of $L=1$~m as in the other two cases.)

\smallskip

In Test case \#2, the horizontal (tangential) component of the velocity is much smaller than the vertical one, except inside the transition region where they become comparable, and the vertical (normal) component of the velocity keeps the same order of magnitude in the whole domain $\Omega_\varepsilon$ where it behaves like $\varepsilon^2$ as discussed earlier in this section. The pressure $p_\varepsilon$ increases significantly inside the porous medium and Fig.~\ref{fig:nff_stokes_circle} shows that this occurs in a continuous way inside a thin region around the upper row of solid inclusions where the pressure also displays the expected constant profile (see Fig.~\ref{fig:nff_sections_epsilon}) except for the oscillations due to the presence of the obstacles.

\medskip

\paragraph{Test case \#3: Oblique forced filtration}\mbox{}\\ 
In this test, the computational domain is $\Omega_\varepsilon = (-0.5,0.5)\times (-0.5,0.5)$~m ($L=1$~m) and we impose $\textbf{f}=(10^{-8},-10^{-7})$~kg/s$^2$, homogeneous Dirichlet boundary conditions on the vertical sides of the domain, homogeneous normal stress $(\mu\nabla\mathbf{u}_\varepsilon - p_\varepsilon\mathbf{I})\,\textbf{n}=\mathbf{0}$~kg/(m\,s$^2$) at $\Gamma_{top} = (-0.5,0.5) \times \{0.5\}$~m and $(\mu\nabla\mathbf{u}_\varepsilon - p_\varepsilon\mathbf{I})\,\mathbf{n}=(0,-10^{-7})$~kg/(m\,s$^2$) at $\Gamma_{bottom}=(-0.5,0.5)\times\{-0.5\}$~m. The characteristic velocity is $U_f\sim 10^{-6}$~m/s and the Reynolds number $Re\sim 1$.
The fluid moves from the top--left corner of the domain towards the bottom--right corner of the domain, but it finds the resistance of the porous media, see Fig.~\ref{fig:test1-2-3}, right.
 
\begin{figure}[h!]
  \hspace*{-6mm}
  \begin{tabular}{ccc}
   $u_{\varepsilon,1}/U_f$ & $u_{\varepsilon,2}/U_f$ & $p_\varepsilon/P_f$ \\
  \includegraphics[width=0.33\textwidth]{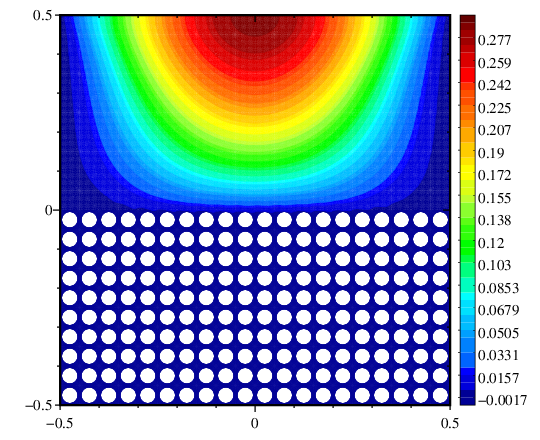} &
  \includegraphics[width=0.33\textwidth]{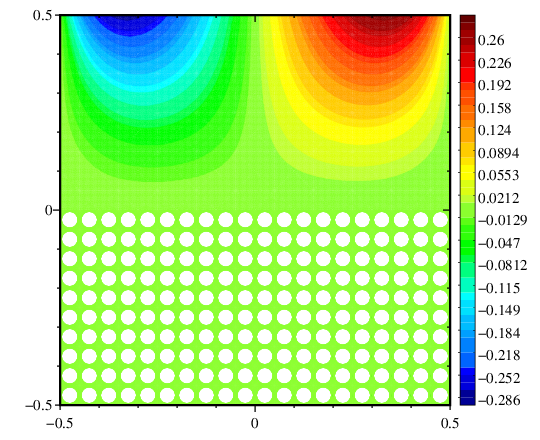} &
  \includegraphics[width=0.33\textwidth]{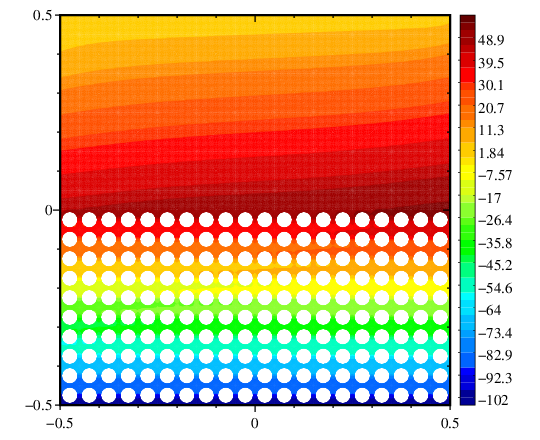}\\
  \includegraphics[width=0.33\textwidth]{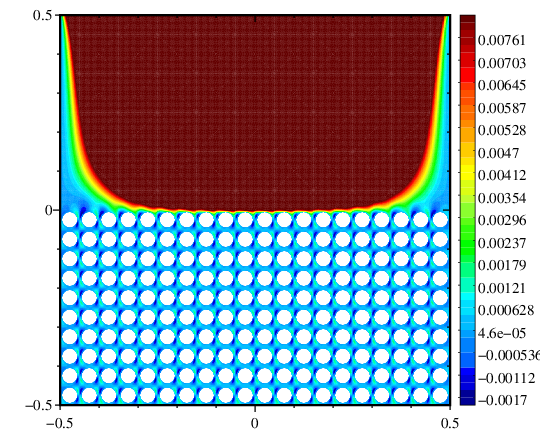} &
  \includegraphics[width=0.33\textwidth]{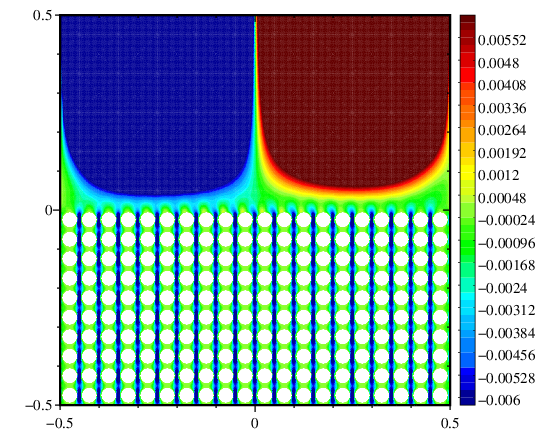} & \\
  \end{tabular}
  \caption{\emph{Test case \#3}. From left to right: horizontal and vertical components of the dimensionless microscale velocity and pressure when $\ell=\frac{1}{20}$ and the circular obstacles have radius $r=0.4\,\ell$. The two components of the velocity in the second row are plotted after changing the colour scale to highlight the small variations of the solution in the porous domain.}
  \label{fig:off_stokes_circle}
\end{figure}

\begin{figure}[h!]
  \centering
  \includegraphics[width=\textwidth]{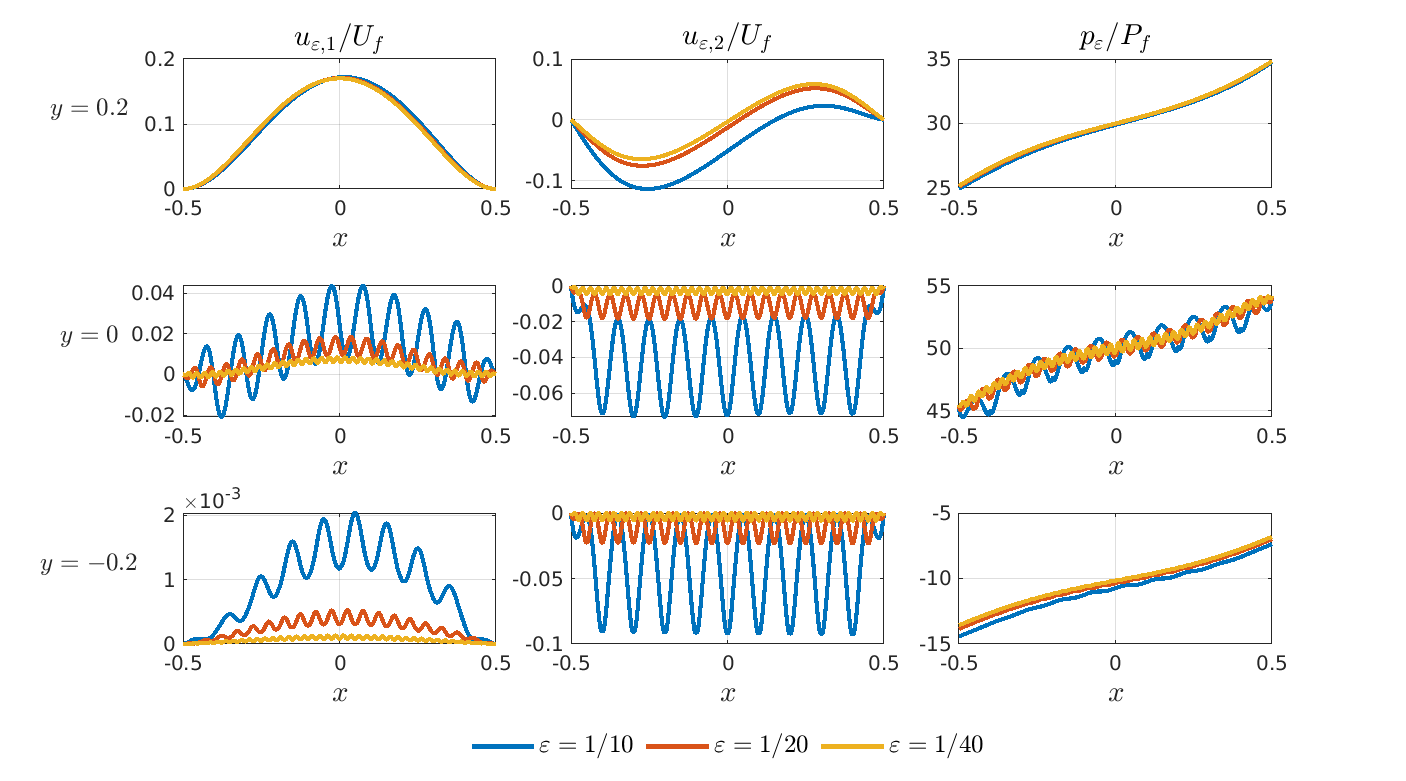}
  \caption{\emph{Test case \#3}. Velocity and pressure profiles of the dimensionless microscale solution at three vertical locations and for three different values of $\varepsilon=\ell$. Circular obstacles with $r=0.4\,\ell$.}
  \label{fig:off_sections_epsilon}
\end{figure}

The velocity and the pressure computed in $\Omega_\varepsilon$ with circular obstacles characterized by $\ell=\frac{1}{20}$ and radius $r=0.4\,\ell$ are shown in Fig.~\ref{fig:off_stokes_circle}, while Fig.~\ref{fig:off_sections_epsilon} reports the dimensionless velocity and pressure profiles for $\ell = \frac{1}{10}, \, \frac{1}{20}, \, \frac{1}{40}$ and $r=0.4\,\ell$ at different vertical locations: $y=0.2$ (in the free-fluid region), $y=0$ (on top of the first row of solid inclusions) and $y=-0.2$ (among the solid obstacles).

\smallskip

Test case \#3 presents an intermediate behaviour with respect to the two previous configurations. Indeed, the results in Figs.~\ref{fig:off_stokes_circle} and \ref{fig:off_sections_epsilon} show that the horizontal component of the velocity decreases throughout the domain $\Omega_\varepsilon$ similarly to its counterpart in Test case \#1. On the other hand, the vertical component of the velocity keeps the same order of magnitude in the domain as it does the corresponding velocity in Test case \#2. The pressure $p_\varepsilon$ keeps the same order of magnitude in $\Omega_\varepsilon$ as it also does in Test case \#1 without the significant variations of Test case \#2 between $\Omega_{\varepsilon f}$ and $\Omega_{\varepsilon p}$.

\bigskip

The physical observations and the numerical results of this section lead us to conclude that the matching conditions \eqref{eq:continuityPressure} and \eqref{eq:continuityVelocity} that impose the continuity of pressure and velocity are physically justified. Also the definition of the global velocity and pressure \eqref{eq:solICDD} is reasonable due to the observation that the Stokes regime still holds in the thin transition region and it is replaced by the Darcy regime only deeper into the porous region.

However, it is crucial that the interfaces $\Gamma_p$ and $\Gamma_f$ are suitably located so that, at the macroscale, the microscopic behaviour of the fluid is correctly captured. This aspect is discussed in the next section.

\subsection{Location of the interfaces $\Gamma_p$ and $\Gamma_f$}\label{sec:Gammaf}

In this section, we provide a practical and physically significant strategy to position the interfaces $\Gamma_p$ and $\Gamma_f$ in such a way that, at the macroscale, the overlapping region between $\Omega_f$ and $\Omega_p$ delimited by the interfaces plays the role of the microscale transition region $\Omega_{\varepsilon t}$ where the sharp variations of velocity and pressure occur.

\smallskip

\begin{figure}[bht]
 \begin{center}
 	\begin{tikzpicture}
        %
        \fill[gray!40] (-0.25, 1.75) circle (4pt);
        \fill[gray!40] ( 0.25, 1.75) circle (4pt);
        \fill[gray!40] ( 0.75, 1.75) circle (4pt);
        \fill[gray!40] ( 1.25, 1.75) circle (4pt);
        \fill[gray!40] ( 1.75, 1.75) circle (4pt);
        \fill[gray!40] ( 2.25, 1.75) circle (4pt);
        \fill[gray!40] ( 2.75, 1.75) circle (4pt);
        \fill[gray!40] ( 3.25, 1.75) circle (4pt);
        \fill[gray!40] ( 3.75, 1.75) circle (4pt);
        \fill[gray!40] ( 4.25, 1.75) circle (4pt);
        \fill[gray!40] (-0.25, 2.25) circle (4pt);
        \fill[gray!40] ( 0.25, 2.25) circle (4pt);
        \fill[gray!40] ( 0.75, 2.25) circle (4pt);
        \fill[gray!40] ( 1.25, 2.25) circle (4pt);
        \fill[gray!40] ( 1.75, 2.25) circle (4pt);
        \fill[gray!40] ( 2.25, 2.25) circle (4pt);
        \fill[gray!40] ( 2.75, 2.25) circle (4pt);
        \fill[gray!40] ( 3.25, 2.25) circle (4pt);
        \fill[gray!40] ( 3.75, 2.25) circle (4pt);
        \fill[gray!40] ( 4.25, 2.25) circle (4pt);
        \draw[black,dotted]  (-0.5, 2.0) -- (4.5, 2.0);
        \draw[black,dotted]  (-0.5, 1.5) -- (4.5, 1.5);
        \draw[black,dotted]  (0.0, 1.35) -- (0.0, 2.5);
        \draw[black,dotted]  (0.5, 1.35) -- (0.5, 2.5);
        \draw[black,dotted]  (1.0, 1.35) -- (1.0, 2.5);
        \draw[black,dotted]  (1.5, 1.35) -- (1.5, 2.5);
        \draw[black,dotted]  (2.0, 1.35) -- (2.0, 2.5);
        \draw[black,dotted]  (2.5, 1.35) -- (2.5, 2.5);
        \draw[black,dotted]  (3.0, 1.35) -- (3.0, 2.5);
        \draw[black,dotted]  (3.5, 1.35) -- (3.5, 2.5);
        \draw[black,dotted]  (4.0, 1.35) -- (4.0, 2.5);
        %
        \draw (-0.5,1.35)--(-0.5,2.7);
        \draw ( 4.5,1.35)--( 4.5,2.7);
        %
        \draw[thick, blue]  (-0.5,2.3)--(4.5,2.3);
        \node[black] at (3.6, 2.0)  {$\Gamma_f$};
        \node[black] at (5.3,2.3) {$y=-\delta^*$};
        \draw[thick, brown] (-0.5,2.5)--(4.5,2.5);
        \node[black] at (0.5, 2.75) {$\Gamma_p$};
        \node[black] at (-1.1,2.5) {$y=0$};
    \end{tikzpicture}
	\end{center}
 \caption{Schematic representation of the location of the interfaces $\Gamma_f$ and $\Gamma_p$. The obstacles are only drawn to provide a reference of where the interfaces are located relative to the microscale domain $\Omega_\varepsilon$.}
 \label{fig:interfaces}
\end{figure}
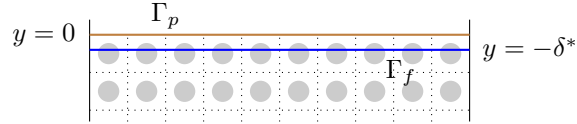

To this aim, first we set the interface $\Gamma_p$ so that, at the microscale, this coincides with the ideal top surface that delimits the solid obstacles of $\Omega_{\varepsilon p}$ as shown in Fig.~\ref{fig:interfaces}, and then letting $\varepsilon \to 0$. 
While from a mathematical point of view, this choice guarantees the periodicity of the microscale representation of the porous medium, from a physical viewpoint, it ensures that the overlapping region is not too low inside the porous region with the risk of getting outside of $\Omega_{\varepsilon t}$.


\smallskip

To identify a correct location for the interface $\Gamma_f$, we consider Test case \#1 (lid--driven cavity problem) as our `model problem' because in this case the variations of the magnitude of the characteristic velocity across $\Omega_\varepsilon$ are much more significant than in the other test cases of Sect.~\ref{sec:microscaleModel}.
The interface $\Gamma_f$ must be accurately placed sufficiently below $\Gamma_p$ to capture the reduction of the magnitude of the velocity up to order $\varepsilon\,U_f$ inside the transition region $\Omega_{\varepsilon t}$. However, it must not be set too low to avoid reaching the bulk region of the porous medium where the Darcy regime characterized by velocity of order $\varepsilon^2\,U_f$ is present. 
Therefore, we look for a rule to characterize the best position of $\Gamma_f$ as a function of known physical parameters (e.g., the porosity and the characteristic microscale length), to guarantee that the ICDD solution approximates the microscale solution with the best possible accuracy.

\smallskip

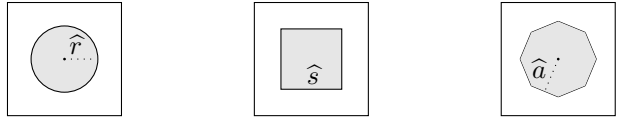
\begin{figure}[bht]
\begin{center}
     \begin{tikzpicture}[scale=0.5]
     \draw (-1.5,-1.5) rectangle (1.5,1.5);
     \fill[draw=black, fill=gray!20] (0,0) circle (25pt);
     \fill[black] (0,0) circle (1pt);
     \draw[dotted] (0,0) -- (0.8,0);
     \node[black] at (0.3, 0.4) {$\widehat{r}$};
     \end{tikzpicture}
     \hspace{1.5cm}
     \begin{tikzpicture}[scale=0.5]
     \draw (-1.5,-1.5) rectangle (1.5,1.5);
     \fill[draw=black, fill=gray!20] (-0.8,-0.8) rectangle (0.8,0.8);
     \node[black] at (0.05, -0.4) {$\widehat{s}$};
     \end{tikzpicture}
     \hspace{1.5cm}
     \begin{tikzpicture}[scale=0.5]
     \tikzmath{\sqrt22 = sqrt(2)/2;}
     \tikzmath{\midx = -sqrt(2)/4;}
     \tikzmath{\midy = (-sqrt(2)/2-1)/2;}
     \draw  (1,0) -- (\sqrt22, \sqrt22);
     \draw  (\sqrt22, \sqrt22)-- (0,1);
     \draw  (0,1)-- (-\sqrt22, \sqrt22);
     \draw  (-\sqrt22, \sqrt22)-- (-1, 0);
     \draw  (-1, 0)-- (-\sqrt22, -\sqrt22);
     \draw  (-\sqrt22, -\sqrt22)-- (0, -1);
     \draw  (0, -1)-- (\sqrt22, -\sqrt22);
     \draw  (\sqrt22, -\sqrt22)-- (1,0) ;
     \fill [fill=gray!20]
           (1,0)
        -- (\sqrt22, \sqrt22)
        -- (0,1)
        -- (-\sqrt22, \sqrt22)
        -- (-1, 0)
        -- (-\sqrt22, -\sqrt22)
        -- (0, -1)
        -- (\sqrt22, -\sqrt22);
     \draw (-1.5,-1.5) rectangle (1.5,1.5);
     \fill[black] (0,0) circle (1pt);
     \draw[dotted] (0,0) -- (\midx,\midy);
     \node[black] at (-0.5, -0.25) {$\widehat{a}$};
    \end{tikzpicture}
\end{center}
\caption{From left to right: reference unit cell $Y$ with circular, square and octagonal obstacles $Y_s$ of radius $\widehat{r}$, side $\widehat{s}$ and apothem $\widehat{a}$, respectively.}
\label{fig:unitCells}
\end{figure}

To be able to compare the microscale solutions of the lid--driven cavity problem with its ICDD counterpart, we first characterize the permeability tensor $\mathbf{K} = K_{ij}$ ($i,j=1,2$) for the Darcy problem \eqref{eq:darcy} following a standard procedure (see, e.g., \cite{Mei:2010}). Namely, let $\mathbf{w}_i = (w_{i1},w_{i2})^T$, for $i=1,2$, be the velocity solution of the auxiliary dimensionless Stokes problem
\begin{equation}\label{eq:reference_cell}
    \begin{array}{rl}
    - \Delta \mathbf{w}_i + \nabla q_i = \mathbf{e}_i & \text{ in } Y_f \\
    \nabla \cdot \mathbf{w}_i = 0\; & \text{ in } Y_f \\
    \mathbf{w}_i & \text{ periodic at } \partial Y \\
    \mathbf{w}_i = \mathbf{0}\; & \text{ at } \partial Y_s \\
    \displaystyle \int_{Y_f} q_i = 0 \,, & 
    \end{array}
\end{equation}
where $\mathbf{e}_i$ is the canonical unit vector in $\mathbb{R}^2$. 
Then, we define the permeability tensor related to the reference cell $Y$ as
\begin{equation*}
    \widehat{K}_{ij} = \int_{Y_f} w_{ij} \, .
\end{equation*}
We consider solid obstacles of three different shapes as illustrated in Fig.~\ref{fig:unitCells} that characterize an isotropic porous medium for which $\widehat{K}_{11}=\widehat{K}_{22}>0$ and $\widehat{K}_{12}=\widehat{K}_{21}=0$, and we set $\widehat K=\widehat K_{11}$.

The computed values of the permeability $\widehat{K}$ for various shapes and sizes of the obstacles $Y_s$ are reported in Table~\ref{tab:obstaclesPorosityOptimaldelta} together with the corresponding porosity $\vartheta = |Y_f|$, being $|Y|=1$. 

Then, let $(\tilde{\mathbf{u}}_\varepsilon, \tilde p_\varepsilon)$ denote the trivial extension of the solution $(\mathbf{u}_\varepsilon, p_\varepsilon)$ of \eqref{eq:stokesGlobal} from $\Omega_\varepsilon$ (with obstacles) to $\Omega$ (without obstacles).

\begin{figure}[h!]
    \centering
    \includegraphics[width=\textwidth]{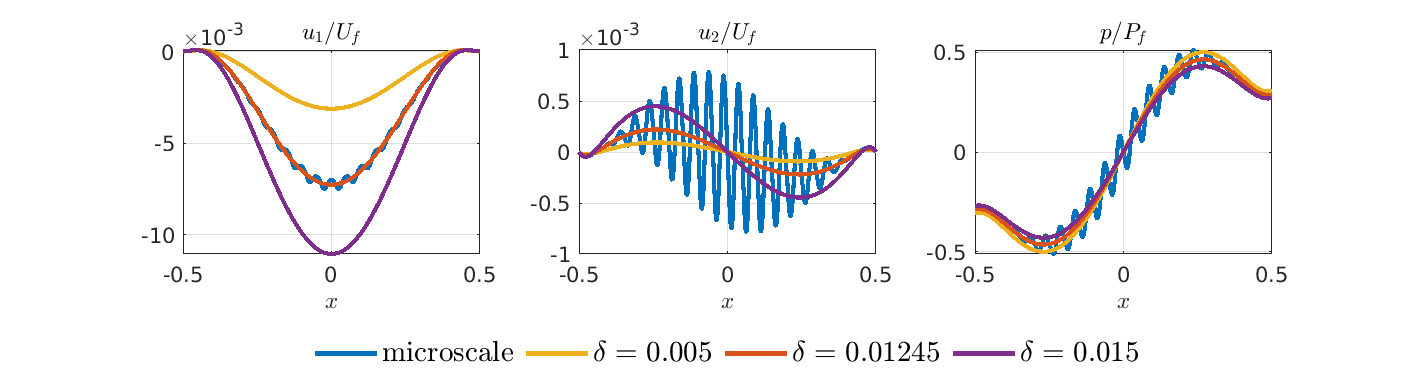}
    \caption{\emph{Test case \# 1}. Comparison between the dimensionless microscale velocity and pressure profiles at $\Gamma_p$ ($y=0$) and the ICDD ones obtained with three values of $\delta$.}
    \label{fig:dnsicdd_y0}
\end{figure}

\begin{figure}[h!]
    \centering
    \includegraphics[width=\textwidth]{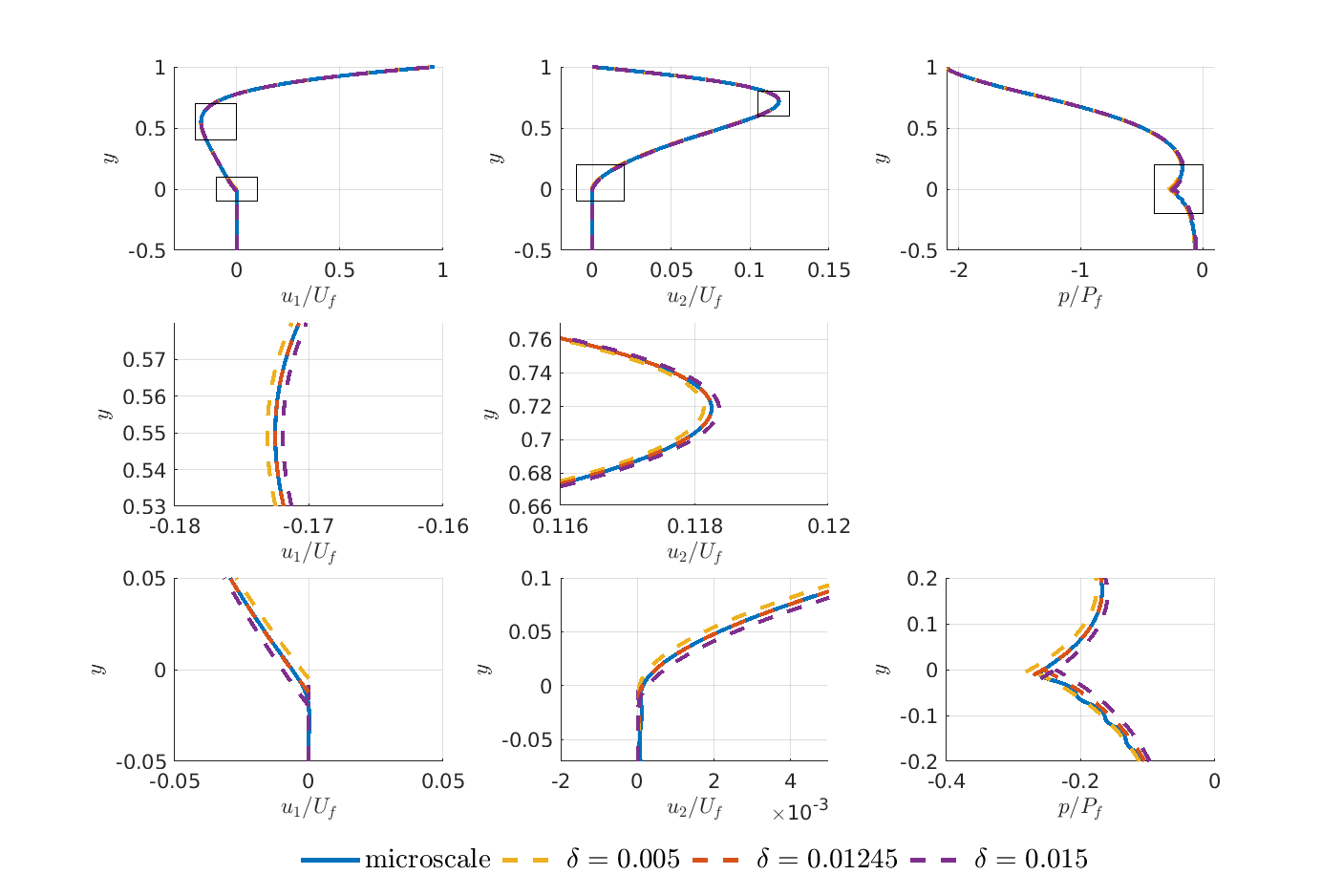}
    \caption{\emph{Test case \#1}. (Top) Dimensionless ICDD velocity and pressure profiles at $x=-0.1$ for three values of $\delta$ and (middle and bottom) zoom of the boxed regions, compared to the microscale solution.}
    \label{fig:dnsicdd_x01}
\end{figure}

\begin{figure}[bht]
    \centering
    \includegraphics[width=\textwidth]{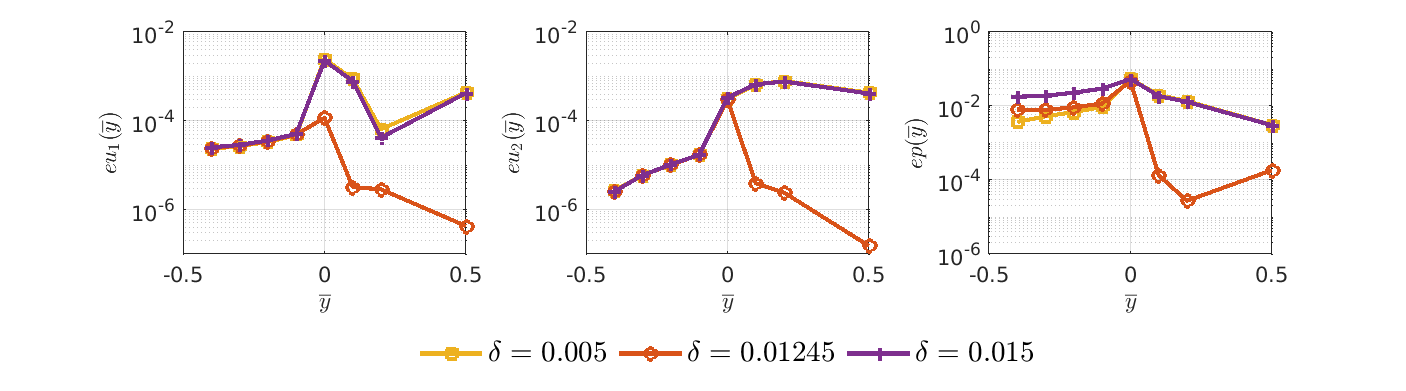}
    \caption{\emph{Test case \# 1.} $L^2$ errors between the ICDD solutions computed with three values of $\delta$ and the microscale solutions at fixed values of $\overline y$. The interface $\Gamma_p$ is at $\overline{y}=0$.}
    \label{fig:dnsicdd_l2norm}
\end{figure}

To understand how the position of the interface $\Gamma_f$ affects the accuracy of the ICDD solution, we perform numerical simulations for Test case \#1 considering the boundary conditions indicated in Sect.~\ref{sec:test1} and circular obstacles. At the microscale, we take $r=\widehat{r}\,\ell$ with $\widehat{r}=0.3$ and $\ell=\frac{1}{20}$ so that the porosity is $\vartheta=0.717$ (see Table~\ref{tab:obstaclesPorosityOptimaldelta}). At the macroscale, we consider the corresponding permeability $K=\ell^2\widehat K$~m$^2$ and we set the interface $\Gamma_f$ at $y=-\delta$ for three different values of $\delta$: $\delta = 0.005, \, 0.01245, \, 0.015$~m. In Fig.~\ref{fig:dnsicdd_y0}, we compare the dimensionless ICDD velocity and pressure profiles with the microscale solution $(\tilde{\mathbf{u}}_\varepsilon, \tilde p_\varepsilon)$ on the interface $\Gamma_p$ ($y=0$), while in Fig.~\ref{fig:dnsicdd_x01} we show the dimensionless profiles along a vertical line at $x=-0.1$. Finally, in Fig.~\ref{fig:dnsicdd_l2norm} we plot the $L^2$ errors
\begin{eqnarray*}
    eu_i(\overline y)&=&\left(
    \int_{-0.5}^{0.5} (\tilde{u}_{\varepsilon,i}(x,\bar{y}) - u_i(x,\bar{y}))^2 \, dx 
    \right)^{\frac{1}{2}}
    \quad \text{for } i=1,2, \text{ and} \\
    ep(\overline y)&=&\left(
    \int_{-0.5}^{0.5} (\tilde{p}_{\varepsilon}(x,\bar{y}) - p(x,\bar{y}))^2 \, dx 
    \right)^{\frac{1}{2}},
\end{eqnarray*}
for $\bar{y}\in\{-0.4,-0.3-0.2,-0.1,0,0.1,0.2,0.5\}$.

The numerical results reported in Figs.~\ref{fig:dnsicdd_y0} and \ref{fig:dnsicdd_x01} show that if $\delta$ is too large, the ICDD Stokes velocity at $y=0$ is underestimated with respect to the microscale one, while the former is overestimated if $\delta$ is too small. This confirms that if the interface $\Gamma_f$ is placed too low inside the porous medium, the ICDD Stokes velocity is incorrectly affected by the Darcy regime $\varepsilon^2\,U_f$, while if $\Gamma_f$ is not low enough, imposing \eqref{eq:continuityVelocity} does not allow to capture the transitional regime of $\Omega_{\varepsilon t}$. Moreover, the errors plotted in Fig.~\ref{fig:dnsicdd_l2norm} indicate that the choice of $\delta$ impacts the behaviour of the Stokes solution more than of the Darcy one so that we aim to determine the best value $\delta^*$ of $\delta$ in order to minimize the difference between the ICDD and the microscale solutions in the free-fluid region. A possible criterion to achieve this is
\begin{equation}\label{eq:criterionDelta0}
    \text{find } \delta^* = \argmin{\delta} ( \| \tilde{\mathbf{u}}_\varepsilon - \mathbf{u}(\delta) \|_{L^2(\Omega_f)} + \|\tilde{p}_\varepsilon - p(\delta) \|_{L^2(\Omega_f)} ) \, ,
\end{equation}
where by $\mathbf{u}(\delta)$ and $p(\delta)$ we denote the dependence of the ICDD solution on $\delta$. 
To avoid numerical artifacts due to the pressure becoming infinite at the top corners of the fluid domain for the lid--driven cavity problem, and also considering that condition \eqref{eq:continuityVelocity} is imposed on the velocity rather than on the pressure, instead of \eqref{eq:criterionDelta0} we adopt the criterion
\begin{equation}
\label{eq:criterionDelta}
    \text{find } \delta^*= \argmin{\delta} \| \tilde{\mathbf{u}}_\varepsilon - \mathbf{u}(\delta) \|_{L^2(\Omega_f^*)}
\end{equation}
where $\Omega_f^* \subset \Omega_f$ is the restricted fluid domain 
\begin{equation}\label{eq:omegaStar}
    \Omega_f^* = \{ (x,y) \in \Omega_f \, : \, -\delta^* \leq y < 0.5 \}\, .
\end{equation}

By an exhaustive method, we compute $\delta^*$ satisfying \eqref{eq:criterionDelta} considering circular obstacles with radius $r=\widehat{r}\,\ell$, square obstacles with side $s=\widehat{s}\,\ell$, and octagonal obstacles of apothem $a=\widehat{a}\,\ell$, for the values of $\widehat{r}$, $\widehat{s}$ and $\widehat{a}$ reported in Table~\ref{tab:obstaclesPorosityOptimaldelta} and for $\ell=\frac{1}{10}$, $\frac{1}{20}$, $\frac{1}{40}$. 

\begin{figure}[bht]
    \centering
    \includegraphics[width=0.5\textwidth]{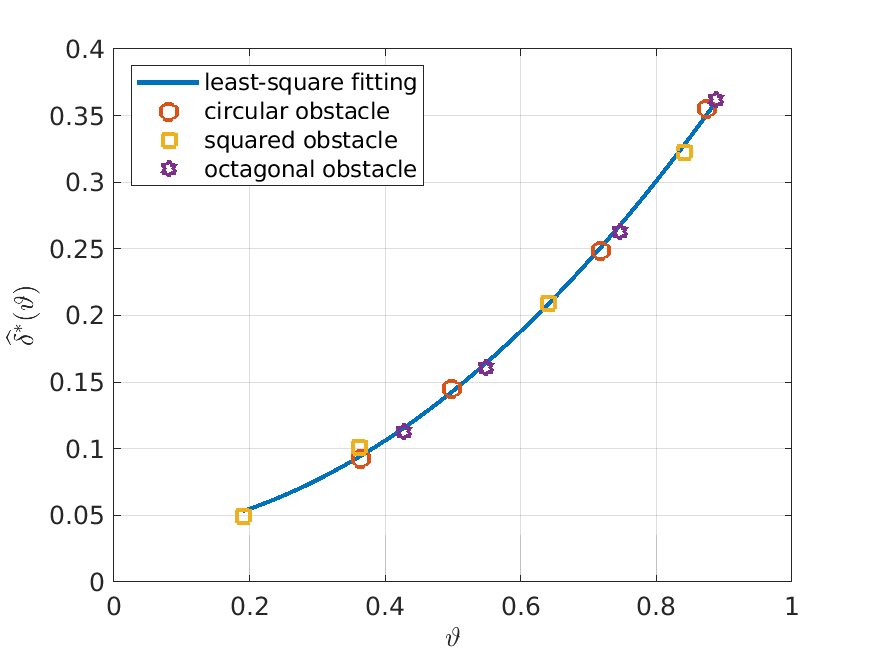}
    \caption{Computed optimal values of $\widehat{\delta}^*$ (symbols) versus the porosity $\vartheta$ for various shapes and sizes of the solid obstacles. The continuous line is the least--square fitting curve \eqref{eq:delta_optim}.}
    \label{fig:deltafitting}
\end{figure}

We observe that the computed optimal values $\delta^*$ are proportional to $\ell$, hence we make them dimensionless by the scaling
\begin{equation*}
    \widehat{\delta}^* = \frac{\delta^*}{\ell}.
\end{equation*}
The resulting optimal values $\widehat{\delta}^*$ are reported in Table~\ref{tab:obstaclesPorosityOptimaldelta}
and plotted in Fig.~\ref{fig:deltafitting} versus the porosity $\vartheta$. We notice that the pairs $(\vartheta,\widehat{\delta}^*)$ describe the same curve independently of the shape of the obstacles and, using least--squares fitting, we obtain the relationship
\begin{equation}\label{eq:delta_optim}
    \widehat{\delta}^*(\vartheta) = 0.3847\,\vartheta^2 + 0.0255\,\vartheta + 0.0344 \, ,
\end{equation}
whose dimensional counterpart becomes
\begin{equation}\label{eq:delta_optim_dimensional}
    \delta^*(\vartheta,\ell) = \ell \, \widehat{\delta}^*(\vartheta)\,.
\end{equation}

The formula \eqref{eq:delta_optim_dimensional} provides a practical way to determine the vertical coordinate 
\begin{equation}\label{eq:Gammaf_position}
y_f=-\delta^*(\vartheta,\ell)
\end{equation}
of the interface $\Gamma_f$ as it only depends on easily computable quantities such as the porosity $\vartheta$ and the characteristic size $\ell$ of the pores of the porous medium. The rule (\ref{eq:Gammaf_position}) is independent of the shape of the obstacles, and it permits to identify a macroscopic overlapping region of thickness proportional to $\ell$, in agreement with theoretical estimates of the dimension of the transition region. Relationship \eqref{eq:delta_optim_dimensional} was obtained considering the lid-driven cavity problem as a model problem because we observed that this was the test case featuring the most complex flow behaviour (see Sect.~\ref{sec:microscaleModel}). However, in Sect.~\ref{sec:validation}, we show that formula \eqref{eq:delta_optim_dimensional} permits to obtain an accurate representation of the flow field in $\Omega$ when $\ell \to 0$ also for the other test cases previously considered.


\begin{table}[bht]
\caption{Computed permeability $\widehat{K}$ and porosity $\vartheta$ for circular, square and octagonal obstacles of different radius $\widehat{r}$, side $\widehat{s}$ and apothem $\widehat{a}$, and optimal $\widehat{\delta}^*$, where $\delta^*$ are computed by solving (\ref{eq:criterionDelta}) on the Test case \#1.}
\label{tab:obstaclesPorosityOptimaldelta}
\begin{center}
\resizebox{\textwidth}{!}{%
\begin{tabular}{lccc|lccc|lccc}
\toprule
\multicolumn{4}{c|}{circular obstacles}& 
\multicolumn{4}{c|}{square obstacles}& 
\multicolumn{4}{c}{octagonal obstacles}\\
\multicolumn{1}{c}{$\widehat{r}$} & $\widehat{K}$ & $\vartheta$ & $\widehat{\delta}^*$ & 
\multicolumn{1}{c}{$\widehat{s}$} & $\widehat{K}$ & $\vartheta$ & $\widehat{\delta}^*$ & 
\multicolumn{1}{c}{$\widehat{a}$} & $\widehat{K}$ & $\vartheta$ & $\widehat{\delta}^*$ \\
\midrule
$0.2$  & $3.295$e$-2$ & $0.874$ & $0.355$ & $0.4$ & $2.358$e$-2$ & $0.840$ & $0.323$ & $0.2$  & $3.561$e$-2$ & $0.887$ & $0.362$\\
$0.3$  & $1.098$e$-2$ & $0.717$ & $0.249$ & $0.6$ & $6.326$e$-3$ & $0.640$ & $0.210$ & $0.3$  & $1.289$e$-2$ & $0.745$ & $0.263$\\
$0.4$  & $1.828$e$-3$ & $0.497$ & $0.145$ & $0.8$ & $7.231$e$-4$ & $0.360$ & $0.102$ & $0.4$  & $2.128$e$-3$ & $0.547$ & $0.161$\\
$0.45$ & $3.173$e$-4$ & $0.364$ & $0.093$ & $0.9$ & $8.651$e$-5$ & $0.190$ & $0.050$ & $0.45$ & $4.378$e$-4$ & $0.427$ & $0.113$\\
\bottomrule
\end{tabular}}
\end{center}
\end{table}

\section{Validation of the ICDD method against direct numerical simulations at the microscale}\label{sec:validation}

In this section, we validate the ICDD model with $\Gamma_p$ and $\Gamma_f$ chosen as described in Sect. \ref{sec:Gammaf} and $\delta^*$ given by (\ref{eq:delta_optim})--(\ref{eq:delta_optim_dimensional}) versus the solutions obtained from direct numerical simulations at the microscale considering the three test cases introduced in Sect.~\ref{sec:microscaleModel} with the same problem data (force and boundary conditions) specified therein. We compare the velocity and the pressure fields computed by these two approaches and we study how the errors behave when $\ell \to 0$ in the microscale setting.

\subsection{Calculation of the errors}
\label{sec:calculationErrors}

To correctly compare the microscale and ICDD velocity fields inside the porous medium, we introduce the affine transformation $F_j^\ell:Y\to Y_j^\ell$ that maps the reference dimensionless cell $Y$ into the generic cell, say $Y_j^\ell$, in $\Omega_\varepsilon$ whose side has length $\ell$. Then, for $i=1,2$, we define the fields ${\bf w}_i^\ell:\Omega_{\varepsilon p}\to {\mathbb R}$ as
\begin{equation*}\label{eq:w_scaled}
\mathbf{w}_i^\ell|_{Y_j^\ell}=\mathbf{w}_i\circ (F_j^\ell)^{-1},
\end{equation*}
where $\mathbf{w}_i$ are the solutions of the auxiliary microscale problem \eqref{eq:reference_cell}. Finally, the function
\begin{equation*}
\tilde{\mathbf{u}}_p = \sum_{i=1}^2 \mathbf{w}_i^\ell u_i
\end{equation*}
is the reconstruction of the velocity in the porous domain combining the ICDD Darcy solution $\mathbf{u}=(u_1,u_2)$ defined in \eqref{eq:solICDD} in $\Omega \setminus \Omega_f$ with the microscale fluctuation $ \mathbf{w}_i^\ell$ (see \cite{Allaire:2010:LN,Allaire:1989:AA}).

We can now define the following errors
\begin{eqnarray}
e^u_{L^2(\Omega_f^*)} = \|\tilde{\mathbf{u}}_\varepsilon -\mathbf{u}\|_{L^2(\Omega_f^*)}, \phantom{e} &
e^p_{L^2(\Omega_f^*)} = \|\tilde p_\varepsilon -p\|_{L^2(\Omega_f^*)}, \label{eq:errors_fluid} \\
e^u_{L^2(\Omega_p^-)} = \|\tilde{\mathbf{u}}_\varepsilon - \tilde{\mathbf{u}}_p 
\|_{L^2(\Omega_p^-)}, &
e^p_{L^2(\Omega_p^-)} = \left\|\tilde p_\varepsilon -p \right\|_{L^2(\Omega_p^-)}, \label{eq:errors_porous}
\end{eqnarray}
where $\Omega_f^*$ is defined in \eqref{eq:omegaStar} while
\begin{equation*}
    \Omega_p^- = \Omega \setminus \Omega_f = \{ (x,y) \in \Omega \, : \, -d \leq y \leq -\delta^* \}\,.
\end{equation*} 



For each of the test cases described in Sect.~\ref{sec:microscaleModel}, we evaluate the errors \eqref{eq:errors_fluid} and \eqref{eq:errors_porous} considering the microscale configurations listed in Table \ref{tab:test_configurations_2}.


\begin{table}
\caption{Configurations for the numerical experiments of Sect.~\ref{sec:validation}.}
\label{tab:test_configurations_2}
\begin{center}
    \begin{tabular}{ccccccc}
    \toprule
    Configuration & Obstacle & Side/Radius & Porosity $\vartheta$ & $\quad\ell\quad$ & $K$ [m$^2$] & $\delta^*$ [m]\\
    \midrule
       &        &              &         & $\frac{1}{10}$ & $7.231$e$-6$ & $9.344$e$-3$\\[2pt]
    C1 & square & $\hat s=0.8$ & $0.360$ & $\frac{1}{20}$ & $1.808$e$-6$ &$4.672$e$-3$\\[2pt]
       &        &              &         & $\frac{1}{40}$ & $4.519$e$-7$ & $2.336$e$-3$\\
    \midrule
       &        &              &         & $\frac{1}{10}$ & $6.326$e$-5$ & $2.083$e$-2$\\[2pt]
    C2 & square & $\hat s=0.6$ & $0.640$ & $\frac{1}{20}$ & $1.582$e$-5$ & $1.041$e$-2$\\[2pt]
       &        &              &         & $\frac{1}{40}$ & $3.954$e$-6$ & $5.207$e$-3$\\
    \midrule
       &        &              &         & $\frac{1}{10}$ & $1.828$e$-5$ & $1.422$e$-2$\\[2pt]
    C3 & circle & $\hat r=0.4$ & $0.497$ & $\frac{1}{20}$ & $4.570$e$-6$ & $7.112$e$-3$\\[2pt]
       &        &              &         & $\frac{1}{40}$ & $1.143$e$-6$ & $3.556$e$-3$\\
    \midrule
       &        &              &         & $\frac{1}{10}$ & $1.098$e$-4$ & $2.506$e$-2$\\[2pt]
    C4 & circle & $\hat r=0.3$ & $0.717$ & $\frac{1}{20}$ & $2.744$e$-5$ & $1.253$e$-2$\\[2pt]
       &        &              &         & $\frac{1}{40}$ & $6.859$e$-6$ & $6.265$e$-3$\\
    \bottomrule
    \end{tabular}
    \end{center}
\end{table}

\paragraph{Numerical comparison between ICDD and Test case \#1 (lid-driven cavity)}\mbox{}

\begin{figure}[h!]
  \hspace*{-6mm}
  \begin{tabular}{ccc}
  $u_{1}/U_f$ & $u_{2}/U_f$ & $p/P_f$ \\
  \includegraphics[width=0.33\textwidth]{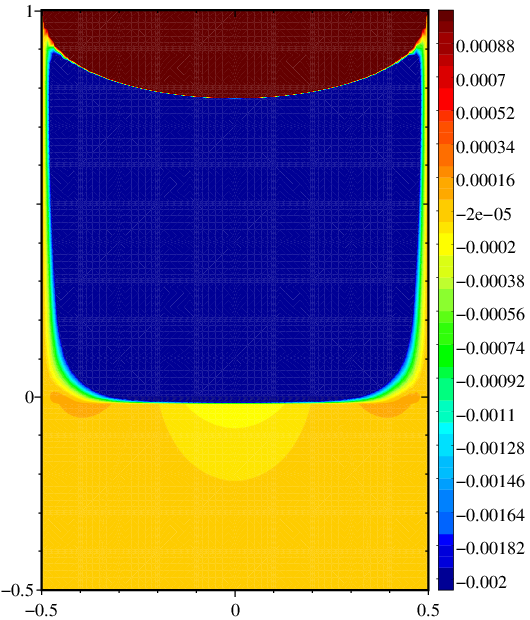} &
  \includegraphics[width=0.33\textwidth]{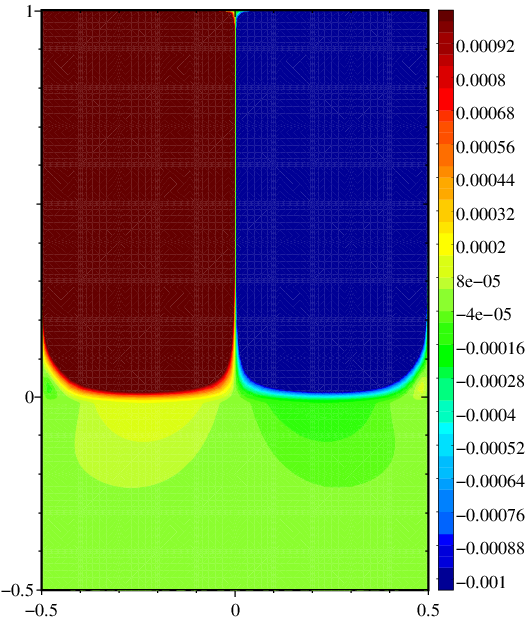} &
  \includegraphics[width=0.33\textwidth]{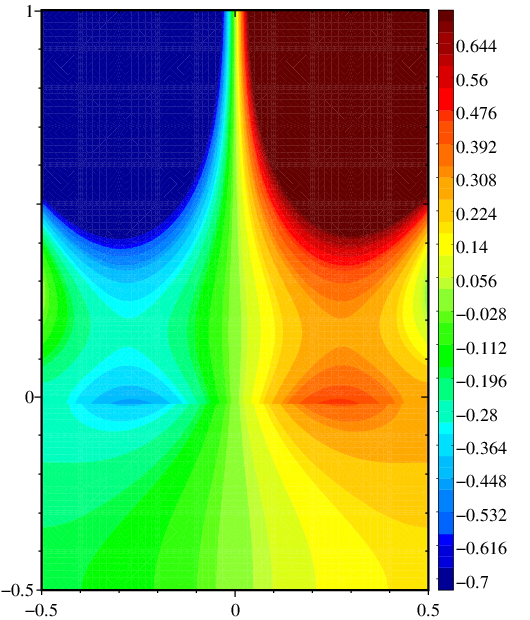}\\
  \end{tabular}
  \caption{\emph{ICDD solution for Test case \#1}. From left to right: horizontal and vertical components of the dimensionless velocity and pressure for circular obstacles with $r=0.2\,\ell$ and $\ell=\frac{1}{20}$.}
  \label{fig:cavity_ICDD}
\end{figure}

We compute the ICDD solution for the lid-driven cavity test, setting the permeability $K=8.237$e$-5$~m$^2$, which corresponds to the case of circular solid inclusions with $r=0.2\,\ell$ and $\ell=\frac{1}{20}$ considered in Sect.~\ref{sec:microscaleModel}, see Fig.~\ref{fig:cavity}. In this case, $\vartheta = 0.874$ and $\delta^* = 1.754$e$-2$~m. The ICDD solution is plotted in Fig.~\ref{fig:cavity_ICDD}, where the values of the colour bar have been chosen to highlight the variations of the solution around the transition region and to easily compare it with the microscale solution of Fig.~\ref{fig:cavity} (bottom). We observe a very good agreement between the macroscopic behaviour of the solutions computed by the two approaches.

\begin{figure}[h!]
  \hspace*{-6mm}
  \includegraphics[width=\textwidth]{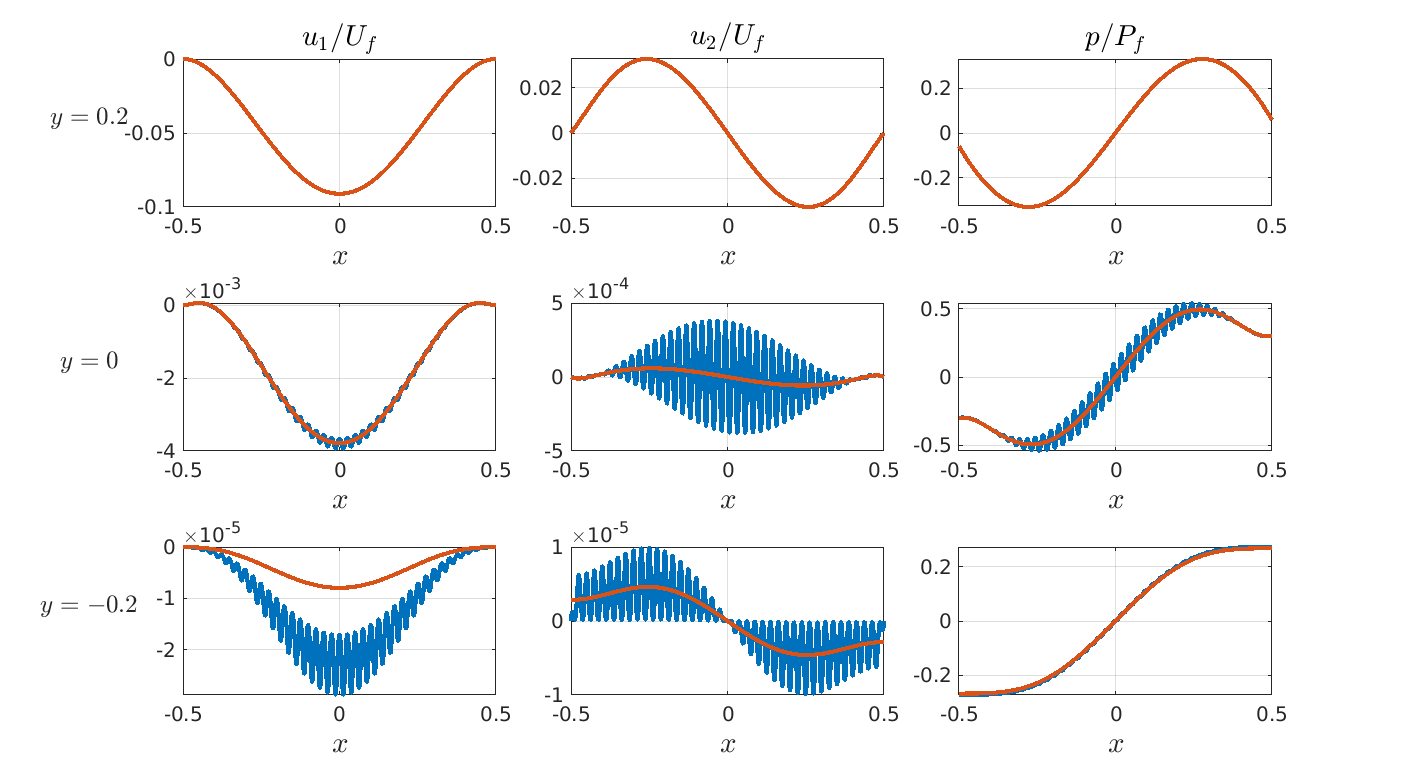} 
  \caption{\emph{Test case \#1}. Dimensionless velocity and pressure profiles of the ICDD solution (red) and reference microscale solution (blue) at three vertical locations for circular obstacles with $r = 0.3\,\ell$ and $\ell=\frac{1}{40}$.}
  \label{fig:sez_ICDD_test1}
\end{figure}

\begin{figure}
\begin{center}
\includegraphics[width=0.9\textwidth]{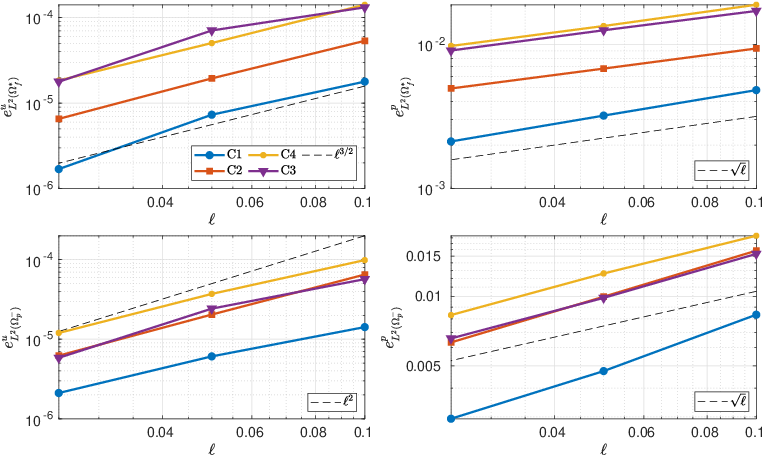}
\end{center}
\caption{\emph{Test case \#1.}  Errors
\eqref{eq:errors_fluid}--\eqref{eq:errors_porous} versus $\ell$ for configurations C1--C4.}
\label{fig:test1_errori_eps}
\end{figure}

Moreover, for configuration C4 (see Table \ref{tab:test_configurations_2}) with $\ell=\frac{1}{40}$, Fig.~\ref{fig:sez_ICDD_test1} shows the traces of the dimensionless microscale (blue) and ICDD (red) solutions on three horizontal lines at $y=0.2$ (in the fluid domain), $y=0$ (on the interface $\Gamma_p$), and $y=-0.2$ (in the porous domain). While the microscale solution presents expected oscillations at the pore scale, the ICDD solution only shows variations at the macroscale, which is correct because the Darcy solution is independent of the microscale. The agreement between the ICDD and the microscale solutions is very good from a qualitative point of view. We can observe that the first component of the ICDD solution inside the porous domain (bottom left panel) underestimates the microscale solution. However, the error estimates in Fig.~\ref{fig:test1_errori_eps} show that, provided that the thickness of the overlapping region is chosen using the least--square fitting formula \eqref{eq:delta_optim_dimensional}, the errors on the velocity  $e^u_{L^2(\Omega_p^-)}$ and $e^u_{L^2(\Omega_f^*)}$ decrease like $\ell^2$ in the porous domain and like $\ell^{3/2}$ in the fluid domain, respectively, while both the errors on the pressure decrease (at least) like $\sqrt{\ell}$.

\begin{figure}[h!]
  \hspace*{-6mm}
  \begin{tabular}{ccc}
  $u_{1}/U_f$ & $u_{2}/U_f$ & $p/P_f$ \\
  \includegraphics[width=0.33\textwidth]{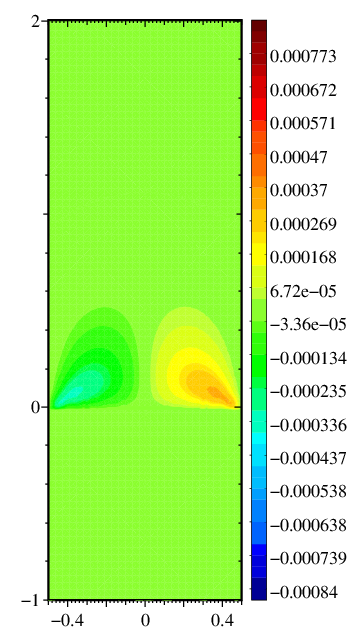} &
  \includegraphics[width=0.33\textwidth]{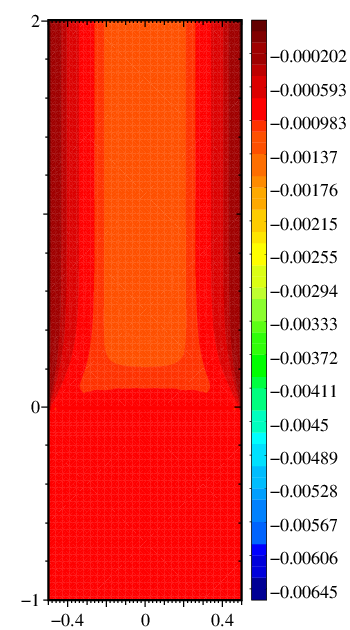} &
  \includegraphics[width=0.33\textwidth]{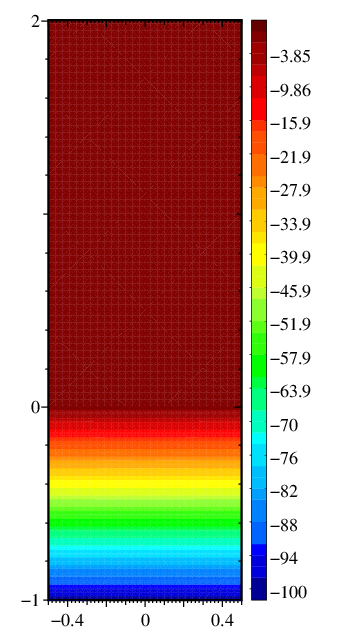}\\
  \end{tabular}
  \caption{\emph{ICDD solution for Test case \#2}. From left to right: horizontal and vertical components of the dimensionless velocity and pressure for circular obstacles with $r=0.4\,\ell$ and $\ell=\frac{1}{20}$.}
  \label{fig:nff_ICDDsolution}
\end{figure}

\begin{figure}[bht]
  \hspace*{-6mm}
       \includegraphics[width=\textwidth]{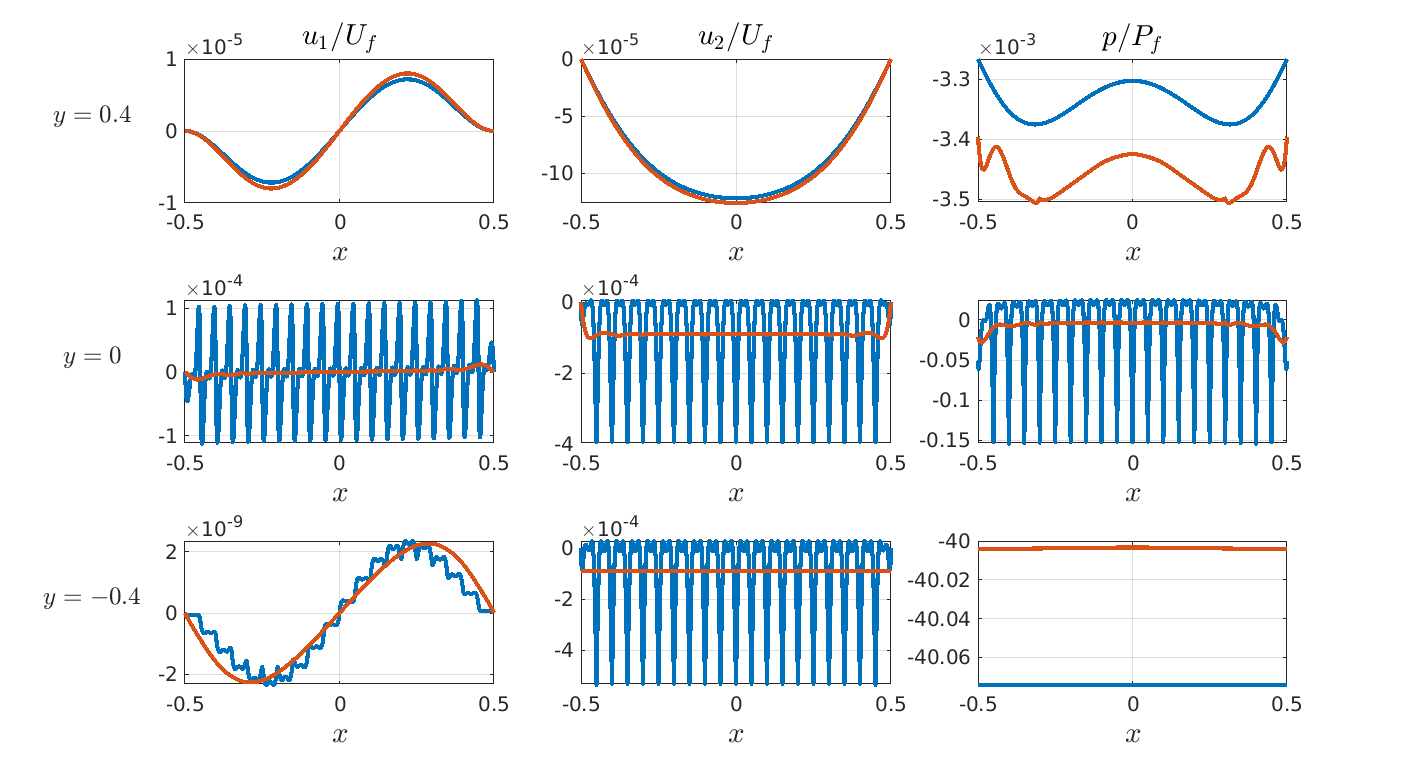}
  \caption{\emph{Test case \#2}. Dimensionless velocity and pressure profiles of the ICDD solution (red) and reference microscale solution (blue) at three vertical locations for square obstacles with $s = 0.8\,\ell$ and $\ell=\frac{1}{40}$.}
  \label{fig:sez_ICDD_test2}
\end{figure}

\paragraph{Numerical comparison between ICDD and Test case \#2 (normal forced filtration)}\mbox{}\\
For this test case, we first compute the ICDD solution using the same setting as in Sect.~\ref{sec:microscaleModel}, see Fig.~\ref{fig:nff_stokes_circle}, which corresponds to configuration C3 (see Table \ref{tab:test_configurations_2}) with $\ell=\frac{1}{20}$.
The numerical results reported in Fig.~\ref{fig:nff_ICDDsolution} show a good agreement with the corresponding microscale solutions of Fig.~\ref{fig:nff_stokes_circle}, and we perform a more quantitative comparison in Fig.~\ref{fig:sez_ICDD_test2}. Here, for configuration C1 with $\ell=\frac{1}{40}$, we plot the dimensionless microscale (blue) and ICDD (red) solutions at $y=0.4$ (in the fluid domain), $y=0$ (on $\Gamma_p$), and $y=-0.4$ (in the porous domain). Similarly to the results obtained for Test case \#1, the ICDD solution only captures variations at the macroscale, but, on average, it accurately represents the behaviour of the solutions observed at the microscale at all three vertical levels.

\begin{figure}
\includegraphics[width=0.9\textwidth]{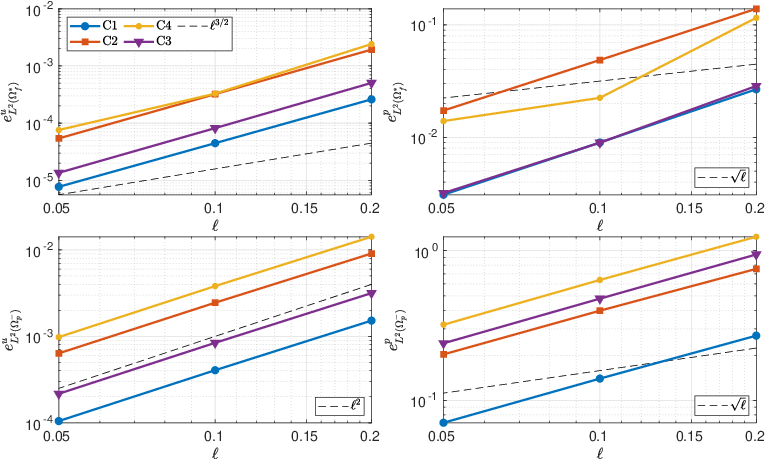}
\caption{\emph{Test case \#2.} Errors \eqref{eq:errors_fluid}--\eqref{eq:errors_porous} versus $\ell$ for configurations C1--C4.}
\label{fig:errori_eps_forced}
\end{figure}

Finally, the numerical study of the errors \eqref{eq:errors_fluid} and \eqref{eq:errors_porous} reported in Fig.~\ref{fig:errori_eps_forced} shows that these decay much faster than in Test case \#1 (refer to Fig.~\ref{fig:test1_errori_eps}). This can be due to the fact that, given the simpler velocity and pressure fields observed for the normal forced filtration (e.g., the whole velocity field just points downwards), the magnitude of the velocity goes to zero very fast when $\ell\to 0$.

\paragraph{Numerical comparison between ICDD and Test case \#3 (oblique forced filtration)}\mbox{}

\begin{figure}[h!]
  \hspace*{-6mm}
  \begin{tabular}{ccc}
  $u_{1}/U_f$ & $u_{2}/U_f$ & $p/P_f$ \\
  \includegraphics[width=0.33\textwidth]{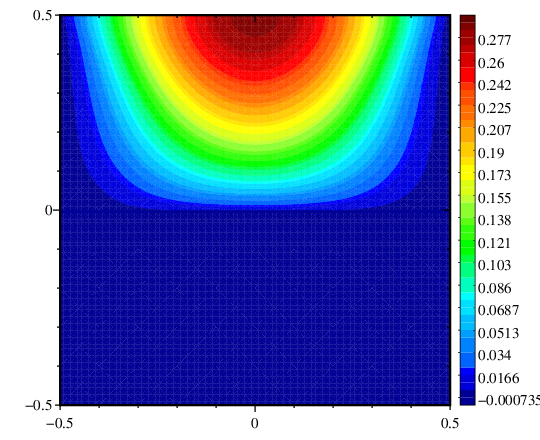} &
  \includegraphics[width=0.33\textwidth]{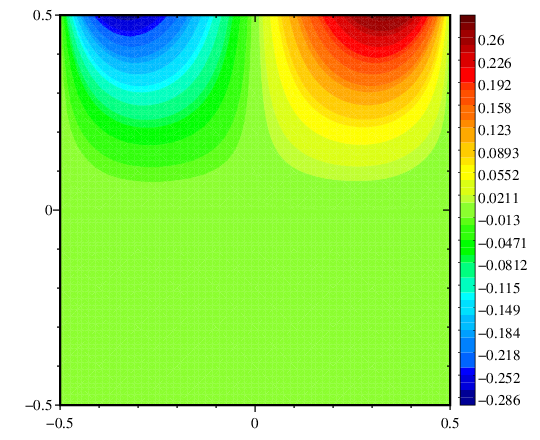} &
  \includegraphics[width=0.33\textwidth]{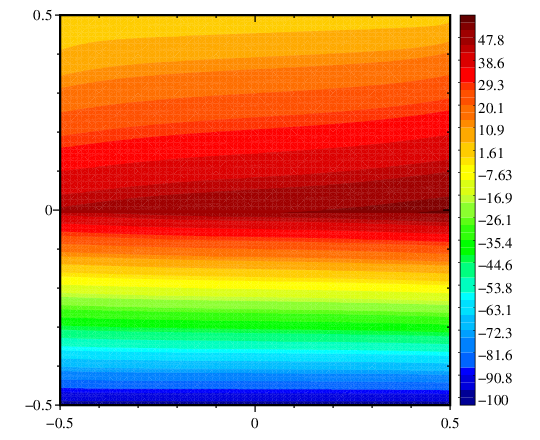}\\
  \end{tabular}
  \caption{\emph{ICDD solution for Test case \#3}. From left to right: horizontal and vertical components of the ICDD velocity and pressure for circular obstacles with $r=0.4\,\ell$ and $\ell=\frac{1}{20}$.}
  \label{fig:off_ICDDsolution}
\end{figure}

Finally, we consider the case of the oblique forced filtration with the same settings used in Sect.~\ref{sec:microscaleModel}, which corresponds to configuration C3 with $\ell = \frac{1}{20}$. Figure~\ref{fig:off_ICDDsolution} shows the two dimensionless components of the velocity and the dimensionless pressure computed by the ICDD method in the whole domain $\Omega$, and we can observe an excellent agreement with the results at the microscale plotted in Fig.~\ref{fig:off_stokes_circle} (top).

\begin{figure}[bht]
  \hspace*{-6mm}
      \includegraphics[width=\textwidth]{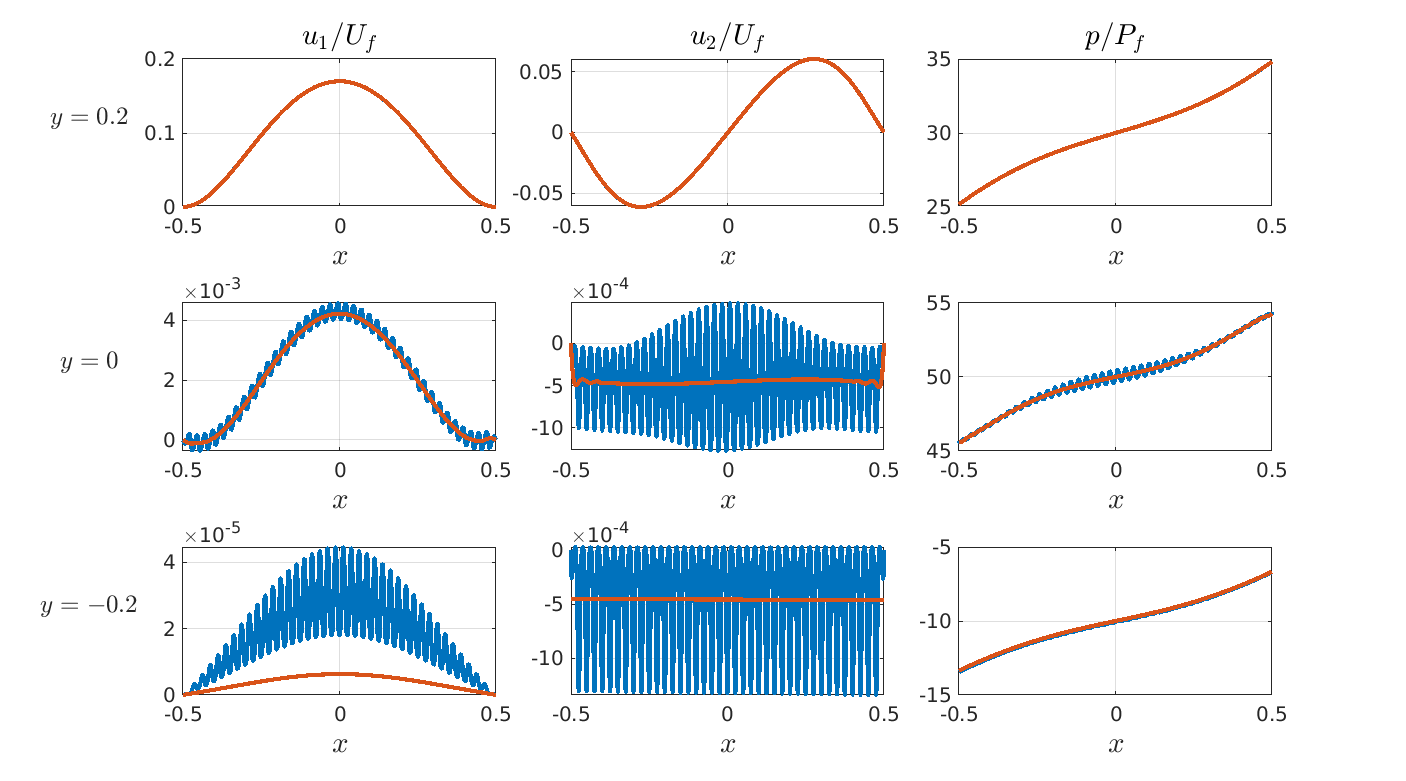} 
  \caption{\emph{Test case \#3}. Dimensionless velocity and pressure profiles of the ICDD solution (red) and reference microscale solution (blue) at three vertical locations for circular obstacles with $r = 0.4\,\ell$ and $\ell=\frac{1}{40}$.}
  \label{fig:sez_ICDD_test3}
\end{figure}

Then, as in the two previous cases, we compare the velocity and pressure profiles obtained by ICDD and by direct numerical simulation at the microscale at the three vertical levels $y=0.2$ (in the fluid domain), $y=0$ (on $\Gamma_p$), and $y=-0.2$ (in the porous domain). For this, we consider configuration C3 with $\ell = \frac{1}{40}$, and we plot the computed solutions in Fig.~\ref{fig:sez_ICDD_test3}. We observe that the ICDD solution correctly represents the macroscopic behaviour of velocity and pressure at the selected levels, with only the horizontal component of the velocity inside the porous medium being slightly underestimated (see bottom left panel).

\begin{figure}
\begin{center}
\includegraphics[width=0.9\textwidth]{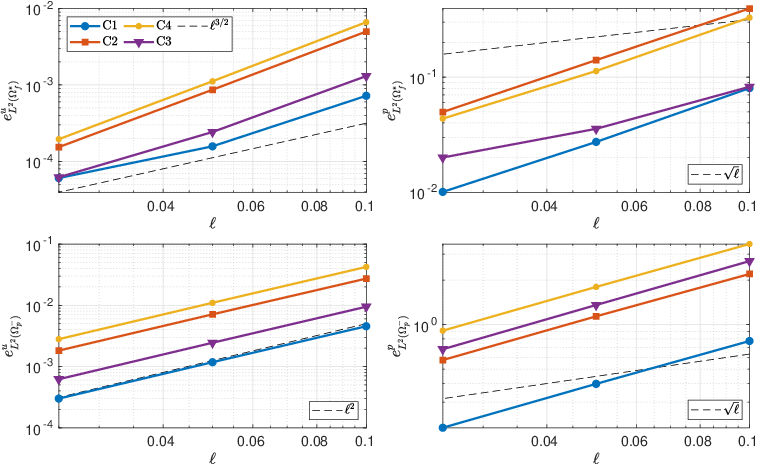}
\end{center}
\caption{\emph{Test case \#3.} Errors \eqref{eq:errors_fluid}--\eqref{eq:errors_porous} versus $\ell$ for configurations C1--C4.}
\label{fig:test3_errori_eps}
\end{figure}

However, the numerical studies of the errors \eqref{eq:errors_fluid}--\eqref{eq:errors_porous} reported in Fig.~\ref{fig:test3_errori_eps} confirm once more that the ICDD solution converges to the microscale solution when the overlap thickness is chosen following the least--square fitting formula \eqref{eq:delta_optim_dimensional}. More precisely, in this case, the errors decrease like $\sqrt{\ell}$ for the pressure (in both the fluid and the porous medium subdomains), while the errors for the velocity decrease like $\ell^2$ inside the porous domain and $\ell^{3/2}$ in the free fluid domain.


\subsection{Error convergence rates}

The numerical study of the errors \eqref{eq:errors_fluid}--\eqref{eq:errors_porous} between the ICDD and the reference microscale solutions carried out in Sect.~\ref{sec:calculationErrors} shows that ICDD solutions converge to the physically correct microscale quantities when $\ell \to 0$. More precisely, the numerical results show that at least approximately the following error convergence rates can be expected:
\begin{equation}\label{eq:convergenceRates}
e^u_{L^2(\Omega_f^*)} \lesssim \ell^{3/2},  \qquad
e^p_{L^2(\Omega_f^*)} \lesssim \sqrt{\ell}, \qquad
e^u_{L^2(\Omega_p^-)} \lesssim \ell^2,      \qquad
e^p_{L^2(\Omega_p^-)} \lesssim \sqrt{\ell},
\end{equation}
with possible higher convergence rates for some particular configurations. This occurs, e.g., in Test cases \#2 and \#3 as it can be seen on the right panels of Figs.~\ref{fig:errori_eps_forced} and \ref{fig:test3_errori_eps}.

\smallskip

Error estimates between a microscale solution and a macroscale (effective) model featuring the Stokes and Darcy problems were previously obtained using homogenization theory in, e.g., \cite{Jager:1996:AnnPisa,Jager:2000:SIAM,Marciniak:2012:MMS,Carraro:2018:NARWA}. (For the sake of clarity, notice that in the cited papers, the characteristic pore size is denoted by $\varepsilon$ instead of by $\ell$ as we do in the present work.)
However, in these works, the effective macroscopic models for the free fluid and the porous medium are coupled in a different way than in the ICDD framework. Indeed, for example, \cite{Marciniak:2012:MMS}, that also considers a flow parallel to the porous medium as in our Test case \#1, macroscopic Stokes and Darcy equations are defined in two non-overlapping subdomains separated by an interface $\Sigma$ whose location coincides with the one of the ICDD interface $\Gamma_p$. Moreover, the coupling between the two local problems occurs only through the Darcy pressure on $\Sigma$ which is set equal to the Stokes pressure plus a correction term proportional to viscous stresses, while the Stokes normal velocity is zero and the tangential velocity satisfies a Beavers-Joseph-type condition on $\Sigma$. In the ICDD framework, two matching conditions \eqref{eq:continuityPressure} and \eqref{eq:continuityVelocity} on both the pressure and the velocity are imposed at two different interfaces, with the position of $\Gamma_f$ depending on $\ell$. Despite these differences, the numerical estimates \eqref{eq:convergenceRates} well align with the theoretical ones of, e.g., \cite{Marciniak:2012:MMS} and references therein.

In fact, \cite{Marciniak:2012:MMS} proved that the interface between the free fluid and the porous medium can be chosen within a layer of pore size thickness and that perturbations of the position of the interface of order $\ell$ only lead to higher-order perturbations on the solution. In the ICDD framework, the interfaces $\Gamma_p$ and $\Gamma_f$ always define a layer of thickness $< \ell$ and the convergence rates observed in all test cases show the robustness of the proposed approach, especially considering that for Test cases \#2 and \#3 we use an optimal value $\delta^*$ to define the interface $\Gamma_f$ that was derived for the case of flow almost parallel to the porous medium. Thus, the ICDD method provides accurate results for isotropic porous media, independently of the flow direction.

Moreover, in the case of fluid flow parallel to the upper surface of the porous medium, \cite{Marciniak:2012:MMS,Jager:2000:SIAM} 
proved that errors in $L^2$ norm between the microscale and macroscale Stokes solution converge like $\sqrt{\ell}$ for the pressure and like $\ell^{3/2}$ for the velocity. We find analogous results (see \eqref{eq:convergenceRates}) for all test cases and not only for the case of the lid-driven cavity (where the flow is almost parallel to the porous medium), and we remark that, differently from the cited works, we impose the (more realistic) coupling condition \eqref{eq:continuityVelocity} on the Stokes velocity at $\Gamma_f$ instead of a no-penetration condition.

Finally, \cite{Carraro:2018:NARWA,Jager:1996:AnnPisa} proved that the microscale pressure in the porous medium strongly converges to the macroscopic Darcy pressure like $\sqrt{\ell}$ with respect to the $L^2$ norm in $\Omega_p$ when a boundary condition on the Cauchy stress is imposed on the top boundary of the periodic microscale porous medium. This boundary condition converges to a macroscale Dirichlet boundary condition for the pressure (see, e.g., \cite{Carraro:2018:NARWA,Fabricius:2017:Cogent}) analogous to the one that we impose on $\Gamma_p$, and we numerically observe the correct convergence rate for the ICDD Darcy pressure. The order of convergence of the $L^2$ norm of the error of the fluid velocity also agrees with the theoretical estimate in \cite{Carraro:2018:NARWA}. 

\begin{figure}
\begin{center}
\includegraphics[width=0.9\textwidth]{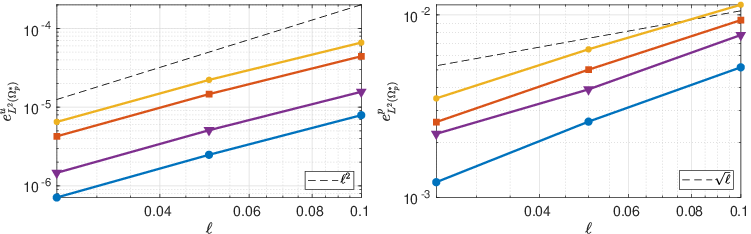}
\end{center}
\caption{\emph{Test case \#1.}  Errors
$e^u_{L^2(\Omega_p^*)}$ and $e^p_{L^2(\Omega_p^*)}$ versus $\ell$ for configurations C1--C4.}
\label{fig:test1_errori_epsbis}
\end{figure}

Notice that for large values of $\ell$, the order of convergence of the ICDD Darcy velocity is slightly less than $\ell^2$. This is particularly evident in Test case \#1 (lid-driven cavity), see Fig.~\ref{fig:test1_errori_eps} (bottom left). This is likely due to the fact that, because of constraints with the computational software, we cannot quantify the boundary layer term of first-order in $\ell$ (see \cite[Theorem 1]{Carraro:2018:NARWA}) which remains non-negligible if $\ell$ is not small enough. In fact, if we replace $\Omega_p^-$ by
\begin{equation}
\Omega_p^* = \{ (x,y) \in \Omega \, : \, -d \leq y \leq -\ell \}\,,
\end{equation}
thus moving away from the boundary layer region close to $\Gamma_f$, the error behaviour significantly improves and we recover the expected convergence rate especially when $\ell$ becomes sufficiently small. This is shown on the left of Fig.~\ref{fig:test1_errori_epsbis}, where we plot the error $e^u_{L^2(\Omega_p^*)}$. We also observe that the errors $e^p_{L^2(\Omega_p^*)}$ computed for the pressure maintain the same behaviour like $\sqrt{\ell}$ as their counterparts in $\Omega_p^-$ (see Fig.~\ref{fig:test1_errori_eps}, bottom right) with only a decrease in magnitude.

\section{Conclusions}

The ICDD method for the coupled Stokes-Darcy problem provides an accurate macroscopic framework for the numerical modelling of the filtration of a laminar incompressible fluid through an isotropic porous medium. Indeed, by carrying out extensive comparisons using reference microscale models, we showed that the macroscale ICDD solution converges to the microscale solution with suitable orders of convergence that depend on the pore scale $\ell$ for both the velocity and the pressure fields, in agreement with available theoretical results based on homogenization theory. In particular, the transitional fluid regime that exists between the free-fluid and the porous-medium regimes is correctly represented in the ICDD approach by the ad-hoc overlapping region defined by a simple formula that only depends on the porosity and the pore scale of the porous medium. This implies that no auxiliary problems must be solved to identify the correct positioning of the interfaces between the fluid and porous medium. The coupling conditions imposed on the interfaces are also independent of any arbitrary parameters, making the setup of the ICDD method very straightforward. Finally, we have shown that the implementation of ICDD is completely non-intrusive so that it can be easily adapted to reusing existing computational software for incompressible porous media flows.

\bigskip

\emph{Acknowledgements.} The first author acknowledges funding through the EPSRC grant EP/V027603/1 and partial support from the QJMAM Fund for Applied Mathematics. The second author was partially supported by PRIN/MUR on grant No.  20227K44ME and by GNCS – INdAM.

\bibliographystyle{plain}
\bibliography{references_icdd}

\begin{thebibliography}{10}

\bibitem{Allaire:2010:LN}
G.~Allaire.
\newblock Homogenization in porous media -- lecture 2 at {CEA}-{EDF}-{INRIA}
  {S}chool on {Homogenization}.
\newblock http://www.cmap.polytechnique.fr/~allaire/homog/, 13-16 December
  2010.

\bibitem{Allaire:1989:AA}
G.~Allaire.
\newblock Homogenization of the {S}tokes flow in a connected porous medium.
\newblock {\em Asymptotic Anal.}, 2(3):203--222, 1989.

\bibitem{Angot:2017:PhysRevE}
P.~Angot, B.~Goyeau, and J.A. Ochoa-Tapia.
\newblock Asymptotic modeling of transport phenomena at the interface between a
  fluid and a porous layer: Jump conditions.
\newblock {\em Phys. Rev. E}, 95:063302, 2017.

\bibitem{Angot:2021:AWR}
P.~Angot, B.~Goyeau, and J.A. Ochoa-Tapia.
\newblock A nonlinear asymptotic model for the inertial flow at a fluid-porous
  interface.
\newblock {\em Adv. Water Resour.}, 149:103798, 2021.

\bibitem{Beavers:1967:JFM}
G.S. Beavers and D.D. Joseph.
\newblock Boundary conditions at a naturally permeable wall.
\newblock {\em J. Fluid Mech.}, 30:197--207, 1967.

\bibitem{Bottaro:2019:JFM}
A.~Bottaro.
\newblock Flow over natural or engineered surfaces: an adjoint homogenization
  perspective.
\newblock {\em J. Fluid Mech.}, 877:877 P1--1--91, 2019.

\bibitem{Bottaro:2020:Meccanica}
A.~Bottaro and S.B. Naqvi.
\newblock Effective boundary conditions at a rough wall: a high-order
  homogenization approach.
\newblock {\em Meccanica}, 55:1781--1800, 2020.

\bibitem{Brinkman:1947:ASRA}
H.C. Brinkman.
\newblock A calculation of the viscous force exerted by a flowing fluid on a
  dense swarm of particles.
\newblock {\em Appl. Sci. Res. A}, 1:27--34, 1947.

\bibitem{Brooks:1982:CMAME}
A.N. Brooks and T.J.R. Hughes.
\newblock Streamline {U}pwind/{P}etrov-{G}alerkin formulations for convection
  dominated flows with particular emphasis on the incompressible
  {N}avier-{S}tokes equations.
\newblock {\em Comput. Meth. Appl. Mech. Eng.}, 32:199--259, 1982.

\bibitem{Carraro:2013:JFM}
T.~Carraro, C.~Goll, A.~Marciniak-Czochra, and A.~Mikeli{\'c}.
\newblock Pressure jump interface law for the {S}tokes-{D}arcy coupling:
  confirmation by direct numerical simulations.
\newblock {\em J. Fluid Mech.}, 732:510--536, 2013.

\bibitem{Carraro:2015:CMAME}
T.~Carraro, C.~Goll, A.~Marciniak-Czochra, and A.~Mikeli{\'c}.
\newblock Effective interface conditions for the forced infiltration of a
  viscous fluid into a porous medium using homogenization.
\newblock {\em Comput. Meth. Appl. Mech. Eng.}, 292:195--220, 2015.

\bibitem{Carraro:2018:NARWA}
T.~Carraro, E.~Maru\v{s}i\'c-Paloka, and A.~Mikeli{\'c}.
\newblock Effective pressure boundary condition for the filtration through
  porous medium via homogenization.
\newblock {\em Nonlinear Anal.-Real World Appl.}, 44:149--172, 2018.

\bibitem{Chandesris:2006:IJHMT}
M.~Chandesris and D.~Jamet.
\newblock Boundary conditions at a planar fluid-porous interface for a
  poiseuille flow.
\newblock {\em Int. J. Heat Mass Transf.}, 49:2137--2150, 2006.

\bibitem{Chandesris:2009:TPM}
M.~Chandesris and D.~Jamet.
\newblock Jump conditions and surface-excess quantities at a fluid/porous
  interface: a multi-scale approach.
\newblock {\em Transp. Porous Med.}, 78:419--438, 2009.

\bibitem{Cimolin:2013:ANM}
F.~Cimolin and M.~Discacciati.
\newblock Navier-{S}tokes/{F}orchheimer models for filtration through porous
  media.
\newblock {\em Appl. Numer. Math.}, 72:205--224, 2013.

\bibitem{Darcy:1856}
H.~Darcy.
\newblock {\em Les {F}ontaines {P}ubliques de la {V}ille de {D}ijon}.
\newblock Dalmont, Paris, 1856.

\bibitem{Discacciati:2016:SINUM}
M.~Discacciati, P.~Gervasio, A.~Giacomini, and A.~Quarteroni.
\newblock The interface control domain decomposition method for
  {S}tokes-{D}arcy coupling.
\newblock {\em SIAM J. Numer. Anal.}, 54(2):1039--1068, 2016.

\bibitem{Discacciati:2013:SICON}
M.~Discacciati, P.~Gervasio, and A.~Quarteroni.
\newblock The interface control domain decomposition ({ICDD}) method for
  elliptic problems.
\newblock {\em SIAM J. Control Optim.}, 51(5):3434--3458, 2013.

\bibitem{Discacciati:2013:JCSMD}
M.~Discacciati, P.~Gervasio, and A.~Quarteroni.
\newblock The interface control domain decomposition ({ICDD}) method for the
  {S}tokes problem.
\newblock {\em J. Coupled Syst. Multiscale Dyn.}, 1(5):372--392, 2013.

\bibitem{Discacciati:2014:IJNMF}
M.~Discacciati, P.~Gervasio, and A.~Quarteroni.
\newblock Interface control domain decomposition methods for heterogeneous
  problems.
\newblock {\em Int. J. Numer. Meth. Fluids}, 76:471--496, 2014.

\bibitem{Discacciati:2009:RMC}
M.~Discacciati and A.~Quarteroni.
\newblock Navier-{S}tokes/{D}arcy coupling: modeling, analysis, and numerical
  approximation.
\newblock {\em Rev. Mat. Complut.}, 22(2):315--426, 2009.

\bibitem{Eggenweiler:2022}
E.~Eggenweiler.
\newblock {\em Interface {C}onditions for {A}rbitrary {F}lows in
  {S}tokes-{D}arcy {S}ystems: {D}erivation, {A}nalysis and {V}alidation}.
\newblock PhD thesis, University of Stuttgart, 2022.

\bibitem{Eggenweiler:2021:MMS}
E.~Eggenweiler and I.~Rybak.
\newblock Effective coupling conditions for arbitrary flows in {S}tokes-{D}arcy
  systems.
\newblock {\em Multiscale Model. Simul.}, 19(2):731--757, 2021.

\bibitem{Ene:1975:JM}
H.I. Ene and E.~S{\'a}nchez-Palencia.
\newblock \'{E}quations et ph\'enomen{\`e}s de surface pour l'\'ecoulement dans
  un mod{\`e}le de milieu poreux.
\newblock {\em J. M{\'e}canique}, 14(1):73--108, 1975.

\bibitem{Fabricius:2017:Cogent}
J.~Fabricius, E.~Miroshnikova, and P.~Wall.
\newblock Homogenization of the {S}tokes equation with mixed boundary condition
  in a porous medium.
\newblock {\em Cogent Math.}, 4:1327502, 2017.

\bibitem{Forchheimer:1901:ZVDI}
P.~Forchheimer.
\newblock Wasserbewegung durch {B}oden.
\newblock {\em Z. Ver. Deutsch. Ing.}, 45:1782--1788, 1901.

\bibitem{Franca:1992:CMAME}
L.P. Franca and S.L. Frey.
\newblock Stabilized finite element methods: {II}. {T}he incompressible
  {N}avier-{S}tokes equations.
\newblock {\em Comput. Meth. Appl. Mech. Engrg.}, 99:209--233, 1992.

\bibitem{Gervasio:1998:NMPDE}
P.~Gervasio and F.~Saleri.
\newblock Stabilized spectral element approximation for the {N}avier--{S}tokes
  equations.
\newblock {\em Numer. Meth. Part Differ. Equ.}, 14(1):115--141, 1998.

\bibitem{Hernandez:2022:CES}
R.~Hernandez-Rodriguez, P.~Angot, B.~Goyeau, and J.A. Ochoa-Tapia.
\newblock Momentum transport in the free fluid-porous medium transition layer:
  one-domain approach.
\newblock {\em Chem. Eng. Sci.}, 248:117111, 2022.

\bibitem{Jager:1996:AnnPisa}
W.~J{\"a}ger and A.~Mikeli\'c.
\newblock On the boundary conditions at the contact interface between a porous
  medium and a free fluid.
\newblock {\em Ann. Scuola Norm. Sup. Pisa Cl. Sci.}, 23:403--465, 1996.

\bibitem{Jager:2000:SIAM}
W.~J{\"a}ger and A.~Mikeli\'c.
\newblock On the interface boundary conditions by {B}eavers, {J}oseph and
  {S}affman.
\newblock {\em SIAM J. Appl. Math.}, 60:1111--1127, 2000.

\bibitem{Jager:2001:SISC}
W.~J{\"a}ger, A.~Mikeli\'c, and N.~Neuss.
\newblock Asymptotic analysis of the laminar viscous flow over a porous bed.
\newblock {\em SIAM J. Sci. Comput.}, 22(6):2006--2028, 2001.

\bibitem{Kang:2024:IJMS}
J.~Kang and M.~Wang.
\newblock Brinkman double-layer model for flow at a free-porous interface.
\newblock {\em Int. J. Mech. Sci.}, 263:108770, 2024.

\bibitem{Kundu:2012}
P.K. Kundu, I.M. Cohen, and D.R. Dowling.
\newblock {\em Fluid {M}echanics}.
\newblock Academic Press, 6th edition, 2016.

\bibitem{Layton:2003:SINUM}
W.L. Layton, F.~Schieweck, and I.~Yotov.
\newblock Coupling fluid flow with porous media flow.
\newblock {\em SIAM J. Num. Anal.}, 40:2195--2218, 2003.

\bibitem{Lacis:2017:JFM}
U.~L{$\bar{\text{a}}$}cis and S.~Bagheri.
\newblock A framework for computing effective boundary conditions at the
  interface between free fluid and a porous medium.
\newblock {\em J. Fluid Mech.}, 812:866--889, 2017.

\bibitem{LeBars:2006:JFM}
M.~{Le Bars} and M.~{Grae Worster}.
\newblock Interfacial conditions between a pure fluid and a porous medium:
  implications for binary alloy solidification.
\newblock {\em J. Fluid Mech.}, 550:149--173, 2006.

\bibitem{Levy:1975:IJES}
T.~Levy and E.~S{\'a}nchez-Palencia.
\newblock On boundary conditions for fluid flow in porous media.
\newblock {\em Int. J. Eng. Sci.}, 13:923--940, 1975.

\bibitem{Marciniak:2012:MMS}
A.~Marciniak-Czochra and A.~Mikeli\'c.
\newblock Effective pressure interface law for transport phenomena between an
  unconfined fluid and a porous medium using homogenization.
\newblock {\em Multiscale Model. Simul.}, 10(2):285--305, 2012.

\bibitem{Masud:2002:CMAME}
A.~Masud and T.J.R. Hughes.
\newblock A stabilized mixed finite element method for {D}arcy flow.
\newblock {\em Comput. Meth. Appl. Mech. Eng.}, 191(39-40):4341--4370, 2002.

\bibitem{Mei:2010}
C.C. Mei and B.~Vernescu.
\newblock {\em Homogenization {M}ethods for {M}ultiscale {M}echanics}.
\newblock World Scientific, 2010.

\bibitem{Naqvi:2021:IJMF}
S.B. Naqvi and A.~Bottaro.
\newblock Interfacial conditions between a free-fluid region and a porous
  medium.
\newblock {\em Int. J. Multiph. Flow}, 141:103585, 2021.

\bibitem{Neale:1974:CJCE}
G.~Neale and W.~Nader.
\newblock Practical significance of {B}rinkman's extension of {D}arcy's law:
  coupled parallel flows within a channel and a bounding porous medium.
\newblock {\em Can. J. Chem. Eng.}, 52(4):475--478, 1974.

\bibitem{OchoaTapia:1995:IJHMT}
J.A. Ochoa-Tapia and S.~Whitaker.
\newblock Momentum transfer at the boundary between a porous medium and a
  homogeneous fluid - {I}. {T}heoretical development.
\newblock {\em Int. J. Heat Mass Transf.}, 38(14):2635--2646, 1995.

\bibitem{Quarteroni:1994}
A.~Quarteroni and A.~Valli.
\newblock {\em Numerical {A}pproximation of {P}artial {D}ifferential
  {E}quations}.
\newblock Springer, Berlin, 1994.

\bibitem{Quarteroni:1999}
A.~Quarteroni and A.~Valli.
\newblock {\em Domain {D}ecomposition {M}ethods for {P}artial {D}ifferential
  {E}quations}.
\newblock Oxford University Press, Oxford, 1999.

\bibitem{Ruan:2023}
L.~Ruan and I.~Rybak.
\newblock Stokes--{B}rinkman--{D}arcy models for coupled free-flow and
  porous-medium systems.
\newblock In E.~Franck, J.~Fuhrmann, V.~Michel-Dansac, and L.~Navoret, editors,
  {\em Finite Volumes for Complex Applications X---Volume 1, Elliptic and
  Parabolic Problems}, pages 365--373. Springer Nature Switzerland, 2023.

\bibitem{Saffman:1971:SAM}
P.G. Saffman.
\newblock On the boundary condition at the interface of a porous medium.
\newblock {\em Stud. Appl. Math.}, 1:93--101, 1971.

\bibitem{Strohbeck:2023:TPM}
P.~Strohbeck, E.~Eggenweiler, and I.~Rybak.
\newblock A modification of the {B}eavers-{J}oseph condition for arbitrary
  flows to the fluid–porous interface.
\newblock {\em Transp. Porous Med.}, 147:605--628, 2023.

\bibitem{Toselli:2005}
A.~Toselli and O.~Widlund.
\newblock {\em Domain {D}ecomposition {M}ethods. {A}lgorithms and {T}heory}.
\newblock Springer, Berlin, 2005.

\bibitem{Valdes:2013:AWR}
F.J. Vald\'es-Parada, C.G. Aguilar-Madera, J.A. Ochoa-Tapia, and B.~Goyeau.
\newblock Velocity and stress jump conditions between a porous medium and a
  fluid.
\newblock {\em Adv. Water Resour.}, 62:327--339, 2013.

\bibitem{Vander:2003}
H.A. van~der Vorst.
\newblock {\em Iterative {K}rylov methods for large linear systems}.
\newblock Cambridge University Press, Cambridge, 2003.

\bibitem{Zampogna:2016:JFM}
G.A. Zampogna and A.~Bottaro.
\newblock Fluid flow over and through a regular bundle of rigid fibres.
\newblock {\em J. Fluid Mech.}, 792:5--35, 2016.

\end{thebibliography}

\end{document}